\documentclass[11pt]{article}
\usepackage{graphicx}
\usepackage{amsfonts}
\usepackage{amssymb}
\usepackage{amsmath,amsthm}
\usepackage{a4wide}
\usepackage{color}
\newcommand{\und}{\underline}
\newcommand{\dis}{\displaystyle}

 \newcommand{\nn}{\nonumber}

\numberwithin{equation}{section}
\theoremstyle{plain}
\newtheorem{theorem}{Theorem}[section]
\newtheorem{lemma}[theorem]{Lemma}

\newtheorem{proposition}[theorem]{Proposition}
\theoremstyle{definition}
\newtheorem{definition}{Definition}[section]

\def\La{{\Lambda}}
\def\ga{{\gamma}}

\def\eps{{\epsilon}}

\def\vep{{\varepsilon}}

\def\square{\ifmmode\sqr\else{$\sqr$}\fi}
\def\sqr{\vcenter{
         \hrule height.1mm
         \hbox{\vrule width.1mm height2.2mm\kern2.18mm\vrule width.1mm}
         \hrule height.1mm}}                  % This is a slimmer sqr.

\title{Super-hydrodynamic limit in interacting particle systems}

\author{
Gioia Carinci$^{\textup{{\tiny(a)}}}$,
Anna De Masi $^{\textup{{\tiny(b)}}}$,\\
Cristian Giardin{\`a} $^{\textup{{\tiny(a)}}}$,
Errico Presutti $^{\textup{{\tiny(c)}}}$.\\\\
{\small $^{\textup{(a)}}$ Dipartimento di Dipartimento di Scienze fisiche, informatiche e matematiche,}\\
{\small Universit\`a di Modena e Reggio Emilia, via Campi 213/b, 41125 Modena, Italy}
\\
{\small $^{\textup{(b)}}$ Dipartimento di Ingegneria e Scienze dell'Informazione e Matematica,}\\
{\small Universit\`a di L'Aquila, via Vetoio 1, 67100 L'Aquila, Italy}\\
{\small $^{\textup{(c)}}$ GSSI, viale F. Crispi 7, 67100 L'Aquila, Italy}\\
{\small }\\
}

\date{\today}

\begin{document}

%\title {Two time scales in macroscopic limits: deterministic  versus stochastic  behavior}
%
\maketitle

\begin{abstract}

This paper is a follow-up of the work initiated in \cite{CDGP1},
where it has been investigated the hydrodynamic limit of
symmetric independent random walkers with birth
at the origin and death at the rightmost occupied site.
Here we obtain two further results: first we characterize the
stationary states on the hydrodynamic time scale and show
that they are given by a family of linear macroscopic profiles
whose parameters are determined by the current reservoirs and
the system mass. Then we prove the existence of a super-hyrdrodynamic
time scale, beyond the hydrodynamic one. On this larger
time scale the system mass fluctuates and correspondingly
the  macroscopic profile of the system randomly moves within the family of  linear profiles,
with the randomness of a Brownian motion.

\end{abstract}

\vskip1cm

\section{Introduction}

In this paper we continue the analysis of the stochastic process introduced in  \cite{CDGP1}. This is a  particles process in the interval $\La_\eps:=[0,\eps^{-1}]\cap \mathbb Z$,
$\eps^{-1}$ a positive integer.
The
%space of the Markov process on the
space of particles configurations is $\mathbb N^{\La_\eps}$,
$\xi = (\xi(x))_{x\in \La_\eps}\in \mathbb N^{\La_\eps}$ and the component $\xi(x)\in\mathbb{N}$ is interpreted as
the number of particles at site $x$. The generator of the Markov process is
%denoted by
 \begin{equation}
\label{generatore}
L   = L^0   +  L_b  + L_{a}
\end{equation}
(dependence on $\eps$ is not made explicit).
$L^0$ is the generator
of the independent random walks process with reflecting boundary conditions,
%it is
%defined on the core of local functions by
   \begin{equation}
\label{0.1}
L^0 f(\xi) =  \frac 12 \sum_{x=0}^{\eps^{-1}-1}
 \xi(x) \left(f(\xi^{x,x+1}) -f(\xi)\right) + \xi(x+1) \left(f(\xi^{x+1,x}) -f(\xi)\right)
%L^0_{x,x+1} f(\xi)
    \end{equation}
%    \begin{equation}
%\label{generatore-bulk}
%L^0_{x,x+1} f(\xi) = \xi(x) \left(f(\xi^{x,x+1}) -f(\xi)\right) + \xi(x+1) \left(f(\xi^{x+1,x}) -f(\xi)\right)
%    \end{equation}
where
$\xi^{x,y}$ denotes the
configuration obtained from $\xi$ by removing
one particle from site $x$ and putting it at site $y$.
%, i.e.
%   \begin{equation*}
%\xi^{x,y}(z) = \left\{
%\begin{array}{ll}
%\xi(z)  & \text{if } z \ne x,y,\\
%\xi(z) -1  & \text{if } z = x,\\
%\xi(z)+1  & \text{if } z=y.
%\end{array} \right.
%   \end{equation*}
  The operator $L_{b} $ describes the action of creating a  particle at the origin at rate $\eps j$, $j>0$:
  \begin{equation}
\label{generatore-birth}
L_{b} f(\xi) = j\eps \left(f(\xi^{+}) -f(\xi)\right),\quad \xi^{+}(x) = \xi(x) + \mathbf{1}_{x=0}\;.
   \end{equation}
Instead $L_a$ removes particles:
    \begin{equation}
\label{generatore-death}
L_{a}f(\xi) =  j\eps \left(f(\xi^{-}) -f(\xi)\right),\quad \xi^{-}(x) = \xi(x) - \mathbf{1}_{x=R_{\xi}}
    \end{equation}
namely a particle is taken out from the edge  $R_{\xi}$ of the configuration $\xi$:
   \begin{equation}
  \label{3.3.1.00}
 \text{ $R_{\xi}$ is such that:}\;\;\;  \begin{cases}
  \xi (y) >0 & \text{for $y= R_{\xi}$} \\
  \xi (y) =0 & \text{for $y>R_{\xi}$\;.} \end{cases}
   \end{equation}
$L_{a}f(\xi)=0$ if $R_{\xi}$ does not exist, i.e.\ if $\xi\equiv 0$.
The removal mechanism is therefore of topological nature, since the determination of the rightmost occupied site requires
a knowledge of the entire configuration. Topological interactions appears in field as diverse as crowd dynamics \cite{piccoli}
or swarm dynamics \cite{parisi}.

\medskip
The independent random walkers process $\{\xi^0_t\}$, i.e. the process with generator $L^0$ and reflecting
boundary conditions at $0$ and $\eps^{-1}$, can be thought as the evolution of  an ``isolated'' system.
The invariant measure for this process (when the total number
$n$ of particles is given) is a product of uniform distributions, {i.e. each of the $n$ particles occupy
each of the $\eps^{-1}+1$ sites with probability $1/(\eps^{-1}+1)$. Moreover}
%and
each particle equilibrates on times $\eps^{-2}t$ (convergence
being
exponentially fast in  $\eps^{-2}t$).

The hydrodynamic limit for such an isolated system
describes the behavior of the particles when $\eps\to 0$: the total number of particles is
taken proportional to $\eps^{-1}$,   times are scaled by
$\eps^{-2}$ while space is scaled down by $\eps$ (so that the macroscopic space
is $[0,1]\subset \mathbb R$). It is well known \cite{dp1} that the limit behavior (under suitable conditions
on the initial configuration)
is then given by the linear heat equation on $[0,1]$ with Neumann boundary conditions
  \begin{equation}
  \label{0.3}
  \frac{\partial\rho}{\partial t} =
  \frac 12 \frac{\partial^2 \rho}{\partial r^2},
  \qquad \;\;\;\frac {\partial \rho}{\partial r}\Big|_{0}
  =\frac {\partial \rho}{\partial r}\Big|_{1} = 0
  \end{equation}
whose solution is $\rho_t(r)  = G_t^{\rm neum}*\rho_0(r)= \int G_t^{\rm neum}(r,r')\rho_0(r')dr'$
where  $G_t^{\rm neum}(r,r')$, $r,r' \in [0,1]$, is
   the Green function
of the heat equation \eqref{0.3}
with
Neumann boundary conditions:
\begin{equation}
\label{e4.3abcdef}
G_t^{\rm neum}(r,r') = \sum_k G_t(r,r'_k),\quad
 G_t(r,r') = \frac{e^{-(r-r')^2/2t}}{\sqrt{2\pi t}}
\end{equation}
$r'_k$ being the images of $r'$ under
repeated reflections of the interval $[0,1]$ to its right and left
(see for instance \cite{KS} pag. 97 for details).
The solution of  \eqref{0.3} converges as $t\to \infty$ exponentially
fast to the uniform distribution.  Thus
the hydrodynamic behavior given by \eqref{0.3}
truly  describes the behavior of the particles  not only on times of
order $\eps^{-2}$ (on which  \eqref{0.3} is derived)
but at all times as well: there is only one time scale in the isolated system. We will see that this is
in contrast with the two time scales in the ``open'' system that we study here,
where ``open'' means that the system is in contact with ``the outside'',
i.e. particles can be created and killed. \\

The type of open systems most studied in the literature is that
 with ``density reservoirs''  \cite{derrida}
which impose an average density $\rho_{+}$ and $\rho_-$ at the boundary sites (respectively $0$ and $\eps^{-1}$)
via creation and annihilation of particles at both sides.  By suitably defining
such birth-death processes,
a system of independent walkers reaches a stationary measure which is a product of
{Poisson} distributions
with average density which interpolates linearly the boundary densities $\rho_{\pm}$,
see {\cite{CGGR} for the finite size correction and also \cite{DF}} where the result is proved for a
class of zero range processes.
In this case the hydrodynamic equation reads
 \begin{equation}
  \label{0.4}
 \frac{\partial\rho}{\partial t} =
  \frac 12 \frac{\partial^2 \rho}{\partial r^2},
  \qquad \;\;\;
    \rho(0)= \rho_+,\;\; \rho(1) = \rho_-
  \end{equation}
and the stationary profile is given by the linear
profile in $[0,1]$ which interpolates between
$\rho_{\pm}$. Again, {also in the $t\to \infty$
limit,} there is complete agreement between the
hydrodynamic equations and the particles process. The system has still  only one time scale.

The density reservoirs creates a non-equilibrium state with a current
flowing through the system.
By the continuity equation such macroscopic current is given by  $\dis{-\frac 12 \frac{\partial \rho}{\partial r}}$
and in the stationary state of equation \eqref{0.4} one recovers Fick's law
$$
\dis{ -\frac 12 \frac{\partial \rho}{\partial r} =  \frac{\rho_+-\rho_-}2}\;.
$$
At the microscopic level, the current generated by the density reservoirs
is the difference between the average number of particles crossing a bond $(x,x+1)$
from the left and  the average number of particles crossing it from the right.  Thus
it is equal to
 \[
 \frac{E[\xi(x)]-E[\xi(x+1)]}{2} \approx \eps\; \frac{\rho_+-\rho_-}2
\]
(denoting here by $E$ expectation with respect to the stationary measure and
recalling that a particle jumps from $x$ to $x+1$ and viceversa at rate $1/2$).  Thus
the micro-current is
proportional to $\eps$.\\

Another option to create a non-equilibrium state in an open
system is to consider ``current reservoirs'' (see also {\cite{dptv1,dptv2,dpt,dfp}}).
They are constructed in such a way
to get directly a current $\eps j$ just
by throwing in particles from
the left at rate $\eps j$ and removing them from the right
at same rate, without fixing the densities at the boundaries.
This is obtained by
the action of $L_b$ in \eqref{generatore-birth} and $L_a$ in \eqref{generatore-death},
which is to add from
the left and respectively remove from the right  particles
at rate $\eps j$.  As a result, the ``current reservoirs''
directly impose a current $\eps j$.

To better appreciate the role of current reservoirs in a non-equilibrium context
it is useful to draw a parallelism with the problem of fixing a macroscopic
quantity in equilibrium, for instance the magnetization in the Ising model.
In that case one has two possibilities: either one introduces an external magnetic
field which select a macroscopic state with the desired magnetization
or one can choose from the very beginning to restrict the
statistical average to the microscopic configurations compatible with the
desired magnetization (micro-canonical ensemble).  In a similar manner,
to impose a given current in non-equilibrium system satisfying Fourier law,
we can either fix the densities at the boundary (using density reservoirs) or,
alternatively, restrict the system evolution to those trajectory with a prescribed
current. This is precisely what the current reservoirs do.\\

In \cite{CDGP1} the hydrodynamic limit of a system
of symmetric independent walkers with current reservoirs,
namely the process
with generator \eqref{generatore}, has been studied. The result
established in that paper is the existence and continuity of the
macroscopic profile when the microscopic process is started from a sufficiently nice
initial configuration. The hydrodynamic scaling limit is characterized
as the separating elements of upper and lower barriers
(we give in section \ref{seconda} a brief account of the results in  \cite{CDGP1}).

In the present paper we further investigate the macroscopic
properties of the system. As a first result we compute
the stationary macroscopic profiles  in the hydrodynamic limit.
We prove they are given by linear functions with slope $-2j$.
Since here the boundary densities are not
fixed we are in a situation with infinitely many such profiles.
The one that is selected by the system is dictated by the
total mass, which  is a conserved
quantity on the time scale $\eps^{-2}t$.
However, on a longer time scale over which fluctuations
of the total mass are allowed, there is not anymore a
privileged profile and indeed the system will explore
different profiles.
Fluctuations of the total mass
will occur on a {\em super-hydrodynamic} time scale.
More precisely the super-hydrodynamic scaling is obtained
by taking $\eps\to 0$ when the initial  number of particles is
taken proportional to $\eps^{-1}$,  times are scaled by
$\eps^{-3}$ while space is scaled down by $\eps$.
We prove that in this limit the macroscopic profiles of the system
moves randomly over the linear profiles with slope $-2j$ and the motion
is controlled by the rescaled total mass which performs
a Brownian motion reflected at the origin.

While in this paper we deal with independent random walkers
we conjecture the phenomenon
of the existence of a super-hydrodynamic scale in interacting
particle systems coupled to current reservoirs to be quite
universal. More precisely we claim
the same phenomenon is to be expected for all systems
(exclusion walkers, zero-range process, inclusion walkers, ...)
which in the  hydrodynamic limit scale to the free boundary
problem given by
\begin{equation}
\label{0.7}
\frac{\partial\rho}{\partial t} =
\frac 12 \frac{\partial^2 \rho}{\partial r^2} + j D_0 - j D_{R(t)}
\end{equation}
where $R(t)$ is the macroscopic counterpart of the edge
introduced in \eqref{3.3.1.00}, $D_0$ denotes a Dirac delta
at the origin corresponding to creation of particles,
$D_{R(t)}$ denotes a Dirac delta
at $R(t)$ corresponding to removal of the rightmost particles.
This free boundary problem will be studied in
\cite{CDGP3} (see also \cite{CG,O}).
The two-time scales observed in our system is reminiscent
of what is found in the context of  processes with
a localized schock, see for instance \cite{BBDP}. The peculiar and maybe
surprising aspect of the super-hydrodynamic limit is the fact
that on the time scale $\eps^{-3}t$ the system show persistent
randomness, while on the hydrodynamic scale $\eps^{-2}t$
the system follows a deterministic evolution.

The paper is organized as follows. In section \ref{seconda},
after recalling the concept of partial order in the sense
of mass transport and the construction of barriers introduced
in \cite{CDGP1}, we state our main results: Theorem \ref{thme4}
which states that the hydrodynamic stationary profiles are the
linear ones; Theorem \ref{thme2.3} describing the profiles
that in the course of time are attracted to the linear ones;
Theorem \ref{Teo:Local} dealing with the super-hydrodynamic
limit. In Section 3 we prove Theorem \ref{thme4}:
we need to perform a separate analysis for the case with
a non-trivial edge ($R({\infty}):=R < 1$) and the case
where the support of the stationary linear profile
coincides with $[0,1]$ $(R=1)$.
In Section 4 we prove the remaining results.
The convergence to linear profiles (Theorem \ref{thme2.3})
is obtained by
introducing a coupling between two processes
and showing that the number of discrepancies
vanishes on the hydrodynamic scale;
the evolution of profiles on the hydrodynamic
time scale is proved by
exploiting the convergence of the law of the mass
density to the law of a Brownian motion on $\mathbb R^+$
reflected at the origin.

\vskip1cm
\section{Definitions and main results}
\label{seconda}

We  consider initial configurations that approximate a macroscopic profile  in the following sense. We first define the local empirical averages of a configuration $\xi\in \mathbb N^{\La_\eps}$ and of a profile $\rho\in L^\infty([0,1],\mathbb R_+)$ as follows. Given  any integer $\ell$  and  $x\in [0,\eps^{-1}-\ell+1]$,  the empirical averages are
%in a configuration $\xi$ as
 \begin{equation}
  \label{3.3.0}
 \mathcal A_\ell(x,\xi) := \frac 1{\ell} \sum_{y=x}^{ x+\ell-1} \xi(y)\qquad \text{and} \qquad
  \mathcal A'_\ell(x,\rho)= \frac 1{\vep\ell} \int_{\eps x}^{ \eps(x+\ell)} \rho(r)dr
\end{equation}

\begin{definition} (Assumptions on the initial conditions.)
\label{def:2.1}
We suppose  $\rho_{\rm init}\in C([0,1],\mathbb R_+)$ and, if it exists, we call
$R(0)=\min\{r: \rho_{\rm init}(r')=0 \; \forall r'\in[r,1] \}$, the ``edge'' of $\rho_{\rm init}$. We fix  $b<1$  suitably close to 1 and $a>0$ suitably small,
%for the sake of definiteness we set $b=9/10$
%and $a=1/20$.
we then denote by $\ell$ the integer part of $\eps^{-b}$
and suppose that for any $\eps>0$ the initial configuration $\xi$ verifies
 \begin{equation}
  \label{3.3.1}
\max_{x\in {[0,\eps^{-1}-\ell+1]}}\Big| \mathcal A_\ell(x,\xi)-
  \mathcal A'_\ell(x,\rho_{\rm init})\Big| \le \eps^a\;.
\end{equation}
We suppose moreover that, if $\rho_{\rm init}$ has an edge $R(0)$, then
  \begin{equation}
  \label{3.3.1.00.bis}
  |\eps R_{\xi}-R(0)|\le \eps^a
   \end{equation}
with $R_{\xi}$ defined in \eqref{3.3.1.00}.
We shall denote by $P^{(\eps)}_\xi$ the law of the process with
generator $L$ given in \eqref{generatore} supported at time 0 by a configuration  $\xi$ as above.
\end{definition}

%\vskip.5cm
\subsection*{\em Hydrodynamic limit.}

\noindent
The following Theorem has been proved in \cite{CDGP1}.

\vskip.3cm
 \begin{theorem}[Existence of hydrodynamic limit]
   \label{Teo:Hydro}

Let $\rho_{\rm init}$ and $\xi$  be as in Definition \ref{def:2.1}.
Then there exists a function $\rho_t(r)\ge 0$,
$t\ge 0$, $r\in [0,1]$, equal to $\rho_{\rm init}$
at time $t=0$,
continuous
in $(r,t)$ and such that
for all  $T>0$,  $\zeta>0$ and  $t\in [0,T]$ the following holds
   \begin{equation}
   \label{5.15}
\lim_{\eps \to 0}P^{(\eps)}_{\xi}\Big[\max_{x\in [0,\eps^{-1}]}|\eps F_\eps(x;\xi_{\eps^{-2} t}) - F(\eps x;\rho_{t})| \le \zeta  \Big] = 1
    \end{equation}
where
\begin{equation}
	\label{ee8.1}
F(r;u) = \int_r^1 {u(r')\, dr'}, \qquad F_\eps(x;\xi):=\sum_{y\ge x} \xi(y)\;.
	\end{equation}	
In particular for all smooth $\phi$ and for all $\zeta>0$ one has
$$
\lim_{\eps\to 0}{P^{(\eps)}_\xi}\left[\Big|\eps \sum_{x} \xi_{\eps^{-2}t}(x) \phi(x) - \int_{0}^1 \phi(r)\rho_t(r) dr \Big| \le \zeta\right] = 1\;.
$$
%
%   \begin{equation}
%   \label{eq:Hydro}
%\lim_{\eps \to 0}P^{(\eps)}_\xi\Big[\max_{x\in[0,\eps^{-1}]} |\eps F_\eps(x;\xi_{\eps^{-2}t}) - F(\eps x;\rho_t)| \le \zeta \Big] = 1 %\qquad {\rm as} \quad \eps\to 0
 %   \end{equation}

\end{theorem}

\vskip.3cm

In  \cite{CDGP1} we have also proved that
the limit  profile $\rho_t$ can be identified as
the separating element between barriers, with the barriers defined as
solutions of discrete Stefan problems.
To explain this result, calling $D_0$  the Dirac delta at $0$, we preliminary define the sets
% $U_\delta \subset \mathcal U$
\begin{eqnarray}
\label{e4.3.0.0}
\mathcal U & := & \Big\{u= c D_0+\rho:  c \ge 0,\;\; \rho\in L^\infty([0,1],\mathbb R_+)\Big\} \nonumber\\
\mathcal U_\delta & := & \Big\{u= c D_0+\rho: \int \rho >j\delta,\; c \ge 0,\;\; \rho\in L^\infty([0,1],\mathbb R_+)\Big\}
\end{eqnarray}
and the {\em cut-and-paste} operator
$K^{(\delta)}:  \mathcal U_\delta\to \mathcal U$
\begin{equation}
\label{e4.3.0.1}
K^{(\delta)} u = j\delta D_0+
\mathbf 1_{r\in [0,R^{(\delta)}_u]}u,\qquad R^{(\delta)}_u:\; \int_{R^{(\delta)}_u}^1u(r)dr= j\delta\;.
\end{equation}

%\begin{definition}[Approximate Stefan problems]
\begin{definition} [Barriers]
\label{defin:e2.3}
Given $u\in L^\infty([0,1],\mathbb R_+)$  with $\int u >0$ we define, for all $\delta$ small enough
so that $u\in \mathcal U_{\delta}$, the ``barriers''  {$S_{n\delta}^{(\delta,\pm)}(u)$},  $n\in \mathbb N$,
as follows: we set
%The $(\delta,\pm)$ evolutions are  defined from initial data in
%$\mathcal U_\delta$ first at times $n\delta$,:
%setting
$S^{(\delta,\pm)}_{0}(u)=u$, and, for $n\ge 1$,
\begin{eqnarray}
\label{e2.10}
&&
%\col{u_{n\delta}^{(\delta,-)}}:=
S^{(\delta,-)}_{n\delta}(u)= K^{(\delta)} G_\delta^{\rm neum} *S^{(\delta,-)}_{(n-1)\delta}(u)
\\&&
%\col{u_{n\delta}^{(\delta,+)}}:=
S^{(\delta,+)}_{n\delta}(u)= G_\delta^{\rm neum} * K^{(\delta)} S^{(\delta,{+})}_{(n-1)\delta}(u)
\nn
\end{eqnarray}
where  $G_t^{\rm neum}(r,r')$, $r,r' \in [0,1]$, $t\ge 0$ is
   the Green function
of the linear heat equation on $[0,1]$
with
Neumann boundary conditions.

\end{definition}

\medskip
The functions $S_{n\delta}^{(\delta,\pm)}$ are obtained by alternating
the map $G^{\rm neum}_{\delta}$ (i.e. the heat kernel)
and the cut and paste map $K^{(\delta)}$
(which takes out a mass $j\delta$ from the right and put it back
at the origin,  the macroscopic counterpart of $L_b {+} L_a$).
It can be easily seen that, unlike the true process $(\xi_t)_{t\ge0}$,
the evolutions
$S_{n\delta}^{(\delta,\pm)}$ conserve the total mass,
that $S_{n\delta}^{(\delta,+)}$ maps $\mathcal
U_\delta$ into $L^\infty$ while  $S_{n\delta}^{(\delta,-)}$
has a singular component ($j\delta D_0$)
plus a $L^\infty$ component.

The evolutions $S_{n\delta}^{(\delta,\pm)}$ define barriers  in the sense of the following partial order.

\medskip
\begin{definition}
\label{defin:ee8.1}
 (Partial order).
%Let $\mathcal U=\big\{u = c D_0+ \rho, c \ge 0, \rho\in L^\infty([0,1],\mathbb R_+)\big\}$,
%then
%For any $u \in \mathcal U$ and  $r\in [0,1]$ we set
%
%, any $\eps$, any $\xi \in \mathbb N^{[0,\eps^{-1}]}$, $r\in [0,1]$, $x\in [0,\eps^{-1}]$ we define:
%moreover
For $u$ and $v$ in the set ${\cal U}$ we define
\begin{equation}
%\label{ineq}
\label{e4.3.0.2.1.000}
u\le v \;\;\qquad\text{iff}\;\;\qquad
F(r; u)    \le F(r; v) \;\;\quad\text{for all $r\in [0,1]$\;.}
\end{equation}
where $F(r;\cdot)$ is defined in \eqref{ee8.1}.
\end{definition}
In \cite{CDGP1} we have proved the following Theorem.

\begin{theorem}[Hydrodynamic limit via barriers]
		\label{Teo:Char}
Let $\rho_t$ be the function of Theorem \ref{Teo:Hydro}.
%\col{with $\rho_{\rm init}\in C([0,1],\mathbb R_+)$   as in Definition \ref{def:2.1},}  then
Then $\rho_t$  is the unique separating element between the barriers {$\{S_{n\delta}^{(\delta,-)}(\rho_{\rm init})\}$ and $\{S_{n\delta}^{(\delta,+)}(\rho_{\rm init})\}$,} namely
for any  $t >0 $, any $r\in [0,1]$  and any $n\in \mathbb N$:
%for any $t>0$
	\begin{equation}
\label{e4.3.0.2.1.00}
S_{t}^{(t2^{-n},-)}(\rho_{\rm init}) \le \rho_t\le S_{t}^{(t2^{-n},+)}(\rho_{\rm init})
	\end{equation}
%for all $\delta,\delta',t$ such that $t=k\delta=k'\delta'$.
in the sense of \eqref{e4.3.0.2.1.000}. Furthermore the lower bound is a non  decreasing function of $n$, the upper bound a non  increasing function of $n$
and
%they become equal in the limit $n\to\infty$, i.e.
	\begin{equation}\label{2.14bis}
\lim_{n \to \infty}  \sup_{r \in [0,1]}\big|F(r; S_t^{(t2^{-n},\pm)}(\rho_{\rm init}))-  F(r;  \rho_t)\big|=0
\end{equation}

\end{theorem}

%   \col{In \cite{CDGP1}, Theorem 2.4, we prove that, as long as $\rho(1,t)>0$,
%   \begin{equation}
%   \label{e2.5}
% \rho(r,t) = G^{\rm neum}_t * \rho_{\rm init} (r) + j \int_0^t \{G^{\rm neum}_s (r,0) - G^{\rm neum}_s (r,1)\}ds,\;\; t\in [0,T]
%    \end{equation}
%   On the other hand, in \cite{CDGP3} we prove that $\rho_t$ is the solution of the following free boundary problem
%   \begin{equation}
%  \label{macro.2}
%\frac{\partial u}{\partial t} = \frac 12 \frac{\partial^2 u}{\partial r^2},\quad r\in (0,R(t)),\quad
%\frac{\partial u}{\partial r}\Big|_{r=0,R_t}= -2j,\;\; u(R_t,t)=0
%  \end{equation}
% with $R(t)$ implicitly defined by the condition  $2j \dot{R}(t)=\frac{1}{2}\, u_{rr}(R(t),t)$,
%  as long as $R(t)\le 1$.
%  }

\vskip.7cm

\subsection*{{\em Stationary profiles in the hydrodynamic time scale.}}
Our first result will be a full characterization of the stationary macroscopic states
in the hydrodynamic limit.

\begin{definition} [{Linear profiles}]
   \label{defin0.2}
We denote by $\mathcal M$ ``the manifold'' of density profiles whose elements
are either of the form (i) $\rho(r) = -2j(r-R) \mathbf 1_{r\le R}$, $R\in (0,1)$; or
(ii) $\rho(r) = -2jr +c$, $c\ge 2j$.  They are conveniently parameterized as
$\rho^{(M)}$, {$M\ge 0$},  where $M$ is defined so that:
   \begin{equation}
   \label{e2.6}
\int_0^1 \rho^{(M)}(r) dr= M, \qquad \quad {\rho^{(0)}\equiv 0}
    \end{equation}
    {In particular case (i) corresponds to $M<j$ and case (ii) is found for $M>j$.}

%We shall write $\rho^{(M;j)}$ when we want to underline the dependence on $j$.

\end{definition}

\medskip
%\col{If $M>j$, then $\rho(r,t)=\rho^{(M)}(r)$  satisfies \eqref{e2.5}, while, if $M<j$ and hence $R<1$, then it solves the free boundary problem \eqref{macro.2}, thus}

\vskip.5cm

  \begin{theorem} [Stationary profiles]
   \label{thme4}

If $\rho_{\rm init} \in \mathcal M$ then $\rho_t =\rho_{\rm init}$ for all $t\ge 0$. %Moreover
%if $\int \rho_{\rm init} = M$ then $\rho_t \to \rho^{(M)}$ in the sense that for any $r\in [0,1]$
%   \begin{equation}
%   \label{e2.666}
%\lim_{t\to \infty}F( r,\rho_t)  =  F(r, \rho^{(M)})
%    \end{equation}

   \end{theorem}
   \medskip

%To have a self consistent proof also in the case $M<j$   we shall not rely on
%the above statements on the free boundary problem and prove directly the statement.
%
%We also have direct results about the existence of the edge $R_t$:

%For $M$ small enough $\rho^{(M)}(r)=0$  in $[R,1]$ ($R$ depends on $M$) so that
%Theorem \ref{thme4} tells us that there are limit evolutions with a free boundary.
%Next theorem extends the statement to some non stationary cases:

%\vskip.5cm
%
% \col{\begin{theorem}[Free boundary non-stationary profiles]
%   \label{thmee5}
%{Let $M:= \int \rho_{\rm init}$} and $j'\le j$, thus
%%where $\rho \le \rho'$ is meant in the sense
%%of \eqref{e4.3.0.2.1.000}
%%that
%%$F(r;\rho) \le F(r;\rho')$ for all $r\in [0,1]$.  Then $\rho_t \le \rho^{(M;j')}$, i.e.\
%   \begin{equation}
%   \label{e2.666.1}
%\text{if} \qquad \rho_{\rm init} \le \rho^{(M;j')} \qquad \text{then} \qquad \rho_t \le  \rho^{(M;j')}\quad \text{for all  $t\ge 0$ }
%    \end{equation}
%    in the sense of \eqref{e4.3.0.2.1.000}.
%\end{theorem}}
%
%
%\medskip
%
%  \col{ Theorem \ref{thmee5} shows in particular that if $M<j'$
%then  $\rho_t(r)=0$ for all $r\ge R'=M/j'$ and  all $t\ge 0$, thus $R_t\le R'<1$ for all $t \ge 0$. Indeed
%$F(R', \rho^{(M;j')})=0$ then by \eqref{e2.666.1} also $F(R',\rho_t)=0$,
%this implies $\rho_t(r)=0$ for almost all $r\ge R'$, hence
%the statement follows by the continuity of
%$\rho_t$.}
%
%\vskip.5cm
%

\subsection*{\em Super-hydrodynamic limit.}  Hydrodynamics  describes the behavior of
the system on times $\eps^{-2} t$ in the limit when $\eps\to 0$.  In our case
the limit evolution is given by $\rho_t$ as  obtained  in Theorem \ref{Teo:Hydro}.
Hydrodynamics predicts convergence to equilibrium:

\vskip.5cm

\vskip.5cm

  \begin{theorem} [Convergence to the stationary profiles]
   \label{thme2.3}

If  $\int _0^1\rho_{\rm init}(r)dr = M$ then $\rho_t \to \rho^{(M)}$
in the sense that %for any $r\in [0,1]$
   \begin{equation}
   \label{e2.666}
\lim_{t\to \infty}\sup_{r\in [0,1]} \Big|F( r;\rho_t)  -  F(r; \rho^{(M)})\Big|=0
    \end{equation}

 \end{theorem}
%
%
%  \begin{theorem} [Convergence to the stationary profiles]
%   \label{thme2.3}
%
%If  $\int _0^1\rho_{\rm init}(r)dr = M$ then $\rho_t \to \rho^{(M)}$
%in the sense that %for any $r\in [0,1]$
%   \begin{equation}
%   \label{e2.666}
%\lim_{t\to \infty}\sup_{r\in [0,1]} \Big|F( r;\rho_t)  -  F(r; \rho^{(M)})\Big|=0
%    \end{equation}
%Moreover, for any $\zeta>0$ and $M>0$ there are $t^*>0$ (which diverges as $\zeta\to 0$)
%and $\eps^*>0$ so that for any $\eps\le \eps^*$
%\begin{equation}  \label{e2.6666}
%\sup_{\xi: |\xi|\le M \eps^{-1}} E_\xi^{(\eps)}\Big[ \max_{x\in [0,\eps^{-1}]}\; |\eps F_\eps(x;  \xi_{\eps^{-2}t^*}) - F(\eps x; \rho^{( \eps|\xi|)})| \Big] \le \zeta
%\end{equation}
%where
%$|\xi|=\sum_{x=0}^{\eps^{-1}}\xi(x)$ is the number
%of particles in the configuration $\xi\in \mathbb N^{\Lambda_\eps}$.
%
%
% \end{theorem}

\medskip

As a consequence of \eqref{e2.666} and if
$\xi$ and $\rho_{\rm init}$ are as in {Definition  \ref{def:2.1}} then
for any $\zeta>0$% and $M>0$ there are $t^*>0$ (which diverges as $\zeta\to 0$)
%and $\eps^*>0$ so that for any $\eps\le \eps^*$
\begin{equation}
 \label{e2.6666x}
\lim_{t\to \infty} \lim_{\eps\to 0}P_\xi^{(\eps)}
\Big[ \max_{x\in [0,\eps^{-1}]}\; |\eps F_\eps(x;  \xi_{\eps^{-2}t}) - F(\eps x; \rho^{( M)})| {\ge \zeta\Big] =0}
\end{equation}
where
$ M = F(0;\rho_{\rm init})$. \eqref{e2.6666x} shows
convergence in the hydrodynamic time scale to the invariant
profiles of the limit evolution.  There is here however
an obvious interchange of limits as the convergence to the invariant profile
 is only after taking the limit $\eps\to 0$.
The true long time particle behavior requires instead the
study of the process $\xi_{\eps^{-2}t_\eps}$ where $t_\eps\to \infty$ as $\eps\to 0$.
If in this limit we obtain something different than
\eqref{e2.6666x} then we say that there are other scales
than the hydrodynamical one, that we call super-hydrodynamic.

\medskip

\begin{theorem}
[{Super-hydrodynamic limit}]
 \label{Teo:Local}

 Let $\xi^{(\eps)}$ be a sequence such that $\eps|\xi^{(\eps)}|\to m>0$ as $\eps\to 0$.  Let $t_\eps$
 be an increasing, divergent sequence, then the process $\xi_{\eps^{-2}t_\eps}$ has two regimes:

 \begin{itemize}

 \item Subcritical. Suppose $\eps t_\eps\to 0$, then
 \begin{equation}
 \label{e2.6666xx}
 \lim_{\eps\to 0}P_{\xi^{(\eps)}}^{(\eps)}
\Big[ \max_{x\in [0,\eps^{-1}]}\; |\eps F_\eps(x;  \xi_{\eps^{-2}t_\eps}) - F(\eps x; \rho^{( m)})| \le \zeta\Big] =1
\end{equation}

\item Critical. Let $t_\eps=t\eps^{-1}$ then
\begin{equation}
\label{SuperHydro}
\lim_{\eps\to 0}P_{\xi^{(\eps)}}^{(\eps)}\left[ \max_{x\in [0,\eps^{-1}]}\; |\eps F_\eps(x;  \xi_{\eps^{-3}t}) - F(\eps x; \rho^{(M^{(\eps)}_t)})| \le \zeta \right] = 1
\end{equation}
where $M^{(\eps)}_t$ converges in law as $\eps\to 0$ to
  $B_{jt}$, where  {$(B_t)_{t\ge 0}$}, $B_0=m$, is the Brownian motion on $\mathbb R_+$ with reflections at the origin.

 \end{itemize}
%Let $\rho_{\rm init}$ and $\xi$ be as  in Definition \ref{def:2.1} and $\dis{M^{(\eps)}_t:=  \eps|\xi_{\eps^{-3}t}|}$. Then
%for  any $t>0$,

\end{theorem}

\medskip

Thus  on a first time scale, i.e.\ the subcritical regime,
the process behaves deterministically, it is
attracted by the manifold $\mathcal M$ and equilibrates to one of the
invariant profiles for the limit evolution, the one with
the same mass.  However on longer times of the order $\eps^{-3}t$
it starts moving stochastically on the manifold $\mathcal M$ where it performs a
Brownian motion. The reason is pretty simple because
the total number $|\xi_t|$ of particles at time
$t$ performs a symmetric random walk reflected at the origin:

\medskip

{\begin{theorem}[Distribution of the particles' number]
\label{thme3.6a}
$|\xi_t| = \sum_{x=0}^{\eps^{-1}}\xi_t(x)$ has the law
of a continuous time random walk  on
$\mathbb N$ which jumps with equal probability by $\pm 1$
after an  {exponential} time of  {parameter} $2j\eps$,
the jumps leading to $-1$ being suppressed.
%As a consequence
%if $\xi$ verifies the conditions in Definition \ref{defin0.1},
%then the process $(M^{(\eps)}_t)_{t\ge 0}$ defined by $M^{(\epsilon)}_t := {\epsilon}|\xi_{\eps^{-3}t}|$,
%converges in law as $\eps\to 0$
%to  $(b_{jt})_{t\ge0}$, where  {$(b_t)_{t\ge 0}$} is the brownian motion on $\mathbb R_+$ with reflections at the origin
%which starts from $ {b_0 =\lim_{\epsilon\to 0} \epsilon |\xi| = \int \rho_{\rm {init}}}$.
\end{theorem}}

\vskip.2cm

\vskip1cm
\section{Stationary macroscopic profiles}
%\label{sec:e9}
\label{nona}

In this section we shall study the fixed point of $S_\delta^{(\delta,-)}$ (see Definition \ref{defin:e2.3}) and their limits as $\delta\to 0$.
We will show that the stationary profiles are linear in this limit.

\subsection{The case $R<1$}
We first analyze the case when the total mass in less than
$j$ that yields profiles with support in $[0,R]$ with $R<1$.

\medskip

\begin{theorem}
\label{ethm14}
For any $R\in (0,1)$ and any $\delta>0$ small enough there is a unique, continuous function
$\rho\ge 0$, hereafter called ``stationary profile'',  with support  in $[0,R]$, $R<1$, and such that
    \begin{equation}
     \label{e4.1}
 S_{\delta}^{(\delta,-)}(j\delta D_0+\rho)=j\delta D_0+\rho
   \end{equation}
Moreover $\rho$ is an increasing function of $R$.

\end{theorem}

\medskip

\noindent
{\bf Proof.} By
%\eqref{e4.0}
\eqref{e2.10}
  \[
  S^{(\delta,-)}_\delta (u)= j\delta D_0+
G^{\rm neum}_\delta * u \, \cdot  \mathbf 1_{r\in [0,x]},\quad x= R^{(\delta)}_{G^{\rm neum}_\delta *u}
\]
If $u$ is a fixed point of $S^{(\delta,-)}_\delta$,
i.e. $S^{(\delta,-)}_\delta(u)=u$, then $u = j \delta D_0 + \rho$ with the support
of $\rho=G^{\rm neum}_\delta * u$ being
the interval $[0,x]$.
As we look for solutions with support in $[0,R]$ we must take $x=R$ and thus get for $\rho$ the equation
    \begin{equation}
    \label{e4.2}
\rho(r) = j \delta G_\delta ^{\rm neum}(0,r)  + \int_0^R dr' G_\delta ^{\rm neum}(r',r) \rho(r'),\quad r\in [0,R]
    \end{equation}
The last condition in \eqref{e4.3.0.1} (with $x\to R$) becomes:
   \begin{equation}
    \label{e4.2.1}
\int_R^1 dr' \int_0^R dr G_\delta ^{\rm neum}(r',r)[\rho(r) +j\delta D_0(r)]= j \delta
    \end{equation}
However  \eqref{e4.2.1} is not an extra condition as it is automatically satisfied if
$\rho$ satisfies \eqref{e4.2}:
%indeed
\begin{eqnarray*}
&&\int_R^1 dr' \int_0^R dr G_\delta ^{\rm neum}(r',r)[\rho(r) +j\delta D_0(r)]\\&&=
\int_0^1 dr' \int_0^R dr G_\delta ^{\rm neum}(r',r)[\rho(r) +j\delta D_0(r)] -
\int_0^R dr' \int_0^R dr G_\delta ^{\rm neum}(r',r)[\rho(r) +j\delta D_0(r)]\\&&=
 j\delta+\int_0^R \rho  - \int_0^R \rho
\end{eqnarray*}

The proof of the theorem is then a consequence of the following lemma.
\medskip

\begin{lemma}
Call
   \begin{equation}
    \label{e4.2.2}
g^0_{{\delta}}(r,r')=\mathbf 1_{r,r'\in [0,R]}G_\delta ^{\rm neum}(r,r'),\;\; g^0_{n\delta}= g^0_{\delta}*\cdots*g^0_{\delta},\;\;\text{($n$ times)}
    \end{equation}
Then the series
\begin{equation}
\label{e4.4}
 j\delta \sum_{n\ge 0} g^0_{(n+1)\delta}(0,r)=:\rho(r),\quad r\in [0,R]
\end{equation}
%(obtained by iteration of
%\eqref{e4.2})
is uniformly
convergent in $r$ and $\delta$, so that $\rho$ (defined by \eqref{e4.4}) is the unique solution of  \eqref{e4.2}.
\end{lemma}

\medskip

\noindent
{\bf Proof.}
To prove convergence we observe that there is a positive constant $a$
such that
%, calling $N_\delta$
%the smallest integer for which $\delta N_\delta \ge 1$,
\begin{equation}
\label{e4.5}
 \sup_{r\in [0,R]}   \int g^0_{\delta N_\delta}(r,r') dr' \le 1-a,\quad N_\delta\in \mathbb N:\;\;\delta (N_\delta-1)< 1 \le \delta N_\delta
\end{equation}
($a$ can be taken as the sup of the probability that a Brownian motion on $\mathbb R$ which
starts at $r\in [0,R]$
is in $(R,1)$ at time $\delta N_\delta$).  We have
\[
g^0_{n\delta }(r,r')\le \frac{c}{\sqrt{n\delta}}, \quad \text{for all $n$}
\]
Then
\[
j\delta \sum_{n= 0}^{N_\delta} g^0_{(n+1)\delta}(0,r) \le c'
\]
It  follows from \eqref{e4.5} that
\begin{equation}
\label{e4.5.1}
 g^0_{n\delta }(r,r')  \le (1-a)^{ k-1} \sup_{r''} g^0_{\delta (m+N_\delta)}(r'',r'),\quad n= kN_\delta+m, \: k \ge 1, \: 0\le m<N_\delta
\end{equation}
with $ g^0_{\delta (m+N_\delta)}(r'',r')\le c''$.  Thus
\[
j\delta \sum_{n>N_\delta} g^0_{(n+1)\delta}(0,r) \le j\delta\sum_{k\ge 1}\sum_{m< N_\delta}c'' (1-a)^{ k-1} \le c'''
\]
this proves the Lemma. \qed

\vskip.5cm

Continuity of $\rho$ follows from \eqref{e4.4}. To prove that $\rho$ increases with $R$ we  use the following representation of the
Green function $g^0$: let $I=[a,a']$ be an interval in $[0,R]$, $J=[R,1]$, $I^*$ and $J^*$ the union of the repeated reflections of $I$ and $J$ around $0$ and $1$.  Then
    \begin{equation}
     \label{e4.6.1}
\int_I  g^0_{n\delta }(r,r')dr' = P_r\Big[B_{n\delta} \in I^*, \;\;B_{k\delta}\notin J^*, k\le n\Big]
   \end{equation}
where
$ P_r$ is the law of the  Brownian motion
$(B_s)_{s\ge 0}$ on $\mathbb R$ which starts from $r \in [0,R]$.  The right hand side is clearly increasing with $R$. This concludes the proof of Theorem \ref{ethm14}.
 \qed

\vskip.5cm

\begin{theorem}
\label{ethm15}
Let $\rho^{(\delta,-)}:=j\delta D_0+\rho$, $\rho$ as in Theorem \ref{ethm14}, then
    \begin{equation}
     \label{e4.6}
\lim_{\delta\to 0} \rho^{(\delta,-)} (r)=2j(R-r),\quad r\in [0,R]
   \end{equation}

\end{theorem}

\medskip

\noindent
{\bf Proof.}  The proof is in two steps: in the first one we prove that the series in \eqref{e4.4}
converges as $\delta\to 0$ (it is, approximately, a Riemann sum of an integral) while in the second step we recognize the
limit to be the linear function in \eqref{e4.6}.    We proceed by proving lower and upper bounds which in the limit $\delta \to 0$ will coincide.

\medskip

\noindent
{\em Lower bound}. Let $I=[a,a'] \subset [0,R]$, then by \eqref{e4.6.1}
    \begin{equation}
     \label{e4.7}
\int_a^{a'}  g^0_{t }(r,r')dr' \ge P_r\Big[ B_t \in \{[a,a']\cup[-a',-a]\}, \;\;\sup_{s\le t}|B_s|\le R\Big],\;\;\; t=n\delta
   \end{equation}
Thus denoting by $G^{\rm Dir}_t(r,r')$ the Green function
of the heat equation  $u_t = \frac 12 u_{rr}$ in $[-R,R]$ with Dirichlet boundary conditions $u(\pm R)=0$:
\begin{equation}
\label{e4.8}
\rho^{(\delta,-)} (r) \ge   j\delta \sum_{n\ge 0} \Big(G^{\rm Dir}_{(n+1)\delta}(0,r)+G^{\rm Dir}_{(n+1)\delta}(0,-r)\Big)
\end{equation}
The right hand side is the Riemann sum of the corresponding integral and due to the uniform convergence of the series proved earlier we have
\begin{equation}
\label{e4.9}
\liminf_{\delta\to 0} \rho^{(\delta,-)} (r) \ge   j  \int_0^\infty \Big(G^{\rm Dir}_{t}(0,r)+G^{\rm Dir}_{t}(0,-r)\Big) \, dt
\end{equation}
{Let $v(s,r)$ be the resolvent of the heat equation with Dirichlet boundary conditions in $[-R,R]$, then $v$ verifies the resolvent  equation  $\frac 1 2 v_{rr}+D_0=sv$. Hence the integral $\dis{\int_0^\infty G^{\rm Dir}_t(0,r)}=v(0,r):=v_0(r)$ is the weak  solution of the problem $\dis{\frac 12  v_{rr} +D_0 =0}$, $v(\pm R)=0$, namely  $v_0(r)=R-|r|$, $r\in [-R,R]$. Then, from \eqref{e4.9},
\begin{equation}
\label{e4.10}
\liminf_{\delta\to 0} \rho^{(\delta,-)} (r) \ge   2j v_0(r)= 2j (R-|r|)
\end{equation}}
this proves
that the lower bound agrees with \eqref{e4.6}.

\medskip

\noindent
{\em Upper bound}.
We first observe that there are positive constants $\alpha$ and $\beta$ so that for all $\delta$ small enough:
    \begin{equation}
     \label{e4.12}
 P_0\Big[\sup_{s\le \delta} |B_s|\ge \delta^{1/2-\alpha} \Big] \le e^{-\beta \delta^{-2\alpha}}
   \end{equation}
and get from \eqref{e4.6.1} with $R_\delta:= R+ \delta^{1/2-\alpha}$
    \begin{equation}
     \label{e4.13}
\int_a^{a'}  g^0_{n\delta }(r,r')dr' \le P_r\Big[B_{n\delta} \in \{[a,a']\cup[-a',-a]\}), \;\;|B_s|\le R_\delta, s\le n\delta \Big]
+ n e^{-\beta \delta^{-2\alpha}}
   \end{equation}
We use \eqref{e4.13} for $n \le \delta^{-2}$ and get, recalling \eqref{e4.5.1},
   \begin{equation}
\label{e4.14}
\rho^{(\delta,-)} (r) \le   j\delta \Big(\sum_{n= 0}^{{\delta^{-2}}} [G^{{\rm Dir}, R_\delta}_{(n+1)\delta}(0,r)+G^{{\rm Dir}, R_\delta}_{(n+1)\delta}(0,-r)] + 2\delta^{-{4}} e^{-\beta \delta^{-2\alpha}}\Big) + c(1-a)^{\delta^{-1}} {+j \delta D_0}
\end{equation}
where $G^{{\rm Dir}, R_\delta}$ is the Green function with Dirichlet boundary conditions in $[-R_\delta,R_\delta]$ and $c$ is a suitable constant.  By  letting $\delta\to 0$ we recover the lower bound, we omit the details.  \qed

%\vskip.5cm

\subsection{The case $R=1$}

The analysis so far covers cases where  the limit  profile is a piecewise  linear function with slope $-2j$ in $[0,R]$, $R<1$, and equal to 0
in $[R,1]$. The  mass is therefore ${jR^2}$, hence the analysis does not apply to cases
where the mass is $> {j}$.
As we shall see a posteriori this corresponds to stationary solutions for the $(\delta,-)$ evolution
having  {support of the form $[0,R]$, with} $R=1-A\delta$.
We are going to prove that the analogue of Theorem \ref{ethm14} holds as well when $R=1-A\delta$, $A>0$ and $\delta$ small enough.

Equations \eqref{e4.2}--\eqref{e4.4} hold unchanged but \eqref{e4.5} needs a new proof.  Calling $P_r$ the law of the Brownian motion $B_t$ on $\mathbb R$ which starts from $r{\in [0,R]}$, we have:
\begin{equation}
\label{e4.15}
  \int g^0_{\delta N_\delta}(r,r') dr' = P_r\Big[ B_{n\delta}\notin J^*, \;\;n =1,..,N_\delta\Big]
\end{equation}
where $J=[1-A\delta,1]$ and $J^*$ is the union of all reflections of $J$.

\medskip

\begin{lemma}

There is $a>0$ so that for any $\delta$ small enough
\begin{equation}
\label{e4.15.1}
% \delta\sum_{n:n\delta \le 1}
 \int g^0_{\delta N_\delta}(r,r') dr' \le 1-a,\quad \text{for all $r\in [0,1-A\delta]$}
%  and all $N$: } \;\frac 12 \le N\delta \le 1
\end{equation}

\end{lemma}

\medskip

\noindent
{\bf Proof.} By \eqref{e4.15}
%the largest integer such that $N\delta \le 1$.
%We obviously have
\begin{eqnarray}
\label{e4.15.15}
 \int g^0_{\delta N_\delta}(r,r') dr' &=&  P_r\Big[ B_{n\delta}\notin J^*, \;\;n =1,..,N_\delta\Big] \le  P_r\Big[ B_{n\delta}\notin J, \;\;n =1,..,N_\delta\Big] \nn\\ &=&  P_r[X=0],\quad X:= \sum_{n=1}^{N_\delta} \mathbf 1_{B_{n\delta} \in J }
\end{eqnarray}
Let
\begin{equation}
\label{e4.17}
p_k=P_r\big[ X=k\big],\quad
M_i =\sum_{k\ge 1} p_k k^i,\;\; i=0,1,2
\end{equation}
so that $P_r[X=0] = 1-M_0$. Hence, by \eqref{e4.15.15}, we can take for $a$
in \eqref{e4.15.1} any lower bound for $M_0$.
We are going to show that
\begin{equation}
\label{e4.18}
M_0 \ge  \frac{M_1^2}{2(2M_2+M_1)} \ge  \frac{M_1^2}{6M_2 }
\end{equation}
We have
\[
M_2 \ge \sum_{k\ge k_0} p_k k^2 \ge k_0  \sum_{k\ge k_0} p_k k
\]
We choose $k_0$ to be the smallest integer so that
\[
\frac{M_2}{k_0} \le \frac{M_1}2,\quad \frac{2M_2}{M_1} \le  k_0 \le \frac{2M_2}{M_1}+1
\]
Then
\[
k_0 \sum_{k\le k_0} p_k \ge \sum_{k\le k_0} k p_k =M_1 - \sum_{k>k_0} kp_k \ge M_1 -\frac{M_2}{k_0} \ge \frac {M_1}2
\]
Thus
\[
{M_0 \ge}\sum_{k\le k_0} p_k \ge \frac {M_1}{2k_0} \ge \frac{M_1^2}{2(2M_2+M_1)}
\]
\eqref{e4.18} is thus proved.

We have
\[
M_1 = \sum_{n=1}^{N_\delta}\int_{J}  \frac{e^{-(r-r')^2/(2\delta n)}}{\sqrt{2\pi \delta n}} dr'
\]
thus
\begin{equation*}
%\label{e4.20}
M_1 \ge  \frac{ e^{-c} \delta}{\sqrt{2\pi \delta}} \sum_{n=N_\delta/2}^{N_\delta} n^{-1/2} \ge \frac{ e^{-c} \delta}{\sqrt{2\pi \delta}} \frac{ \sqrt {N_\delta}}{2} \ge C_1
\end{equation*}
where $c$ and $C_1$ are constant independent of $r$ and  $\delta$ (recall that $N \approx \delta^{-1}$).  An analogous proof yields $M_1\le C_2$, $C_2$ a constant independent of $r$ and $\delta$.  Moreover
\[
M_2 = M_1 +\sum_{1\le n_1 < n_2  \le N_\delta}\int_{J}dr' \int_{J}dr'' \;\; \frac{e^{-(r-r')^2/(2\delta n_1)}}{\sqrt{2\pi \delta n_1}}  \;\;\frac{e^{-(r'-r'')^2/(2\delta n_2)}}{\sqrt{2\pi \delta n_2}}
\]
As before we can prove (details are omitted) that $M_2 \le C_3$,  a constant independent of $r$ and  $\delta$.

Since  $P_r[X=0] = 1-M_0$
the above together with \eqref{e4.15.15} proves   the lemma with
\begin{equation*}
%\label{e4.20}
a =  \frac{C_1^2}{6C_3}
\end{equation*}
\qed

\vskip.5cm

After \eqref{e4.5} the proof of Theorem \ref{ethm14} extends unchanged to the present case, so that
the conclusions of  Theorem \ref{ethm14} hold as
well when $R=1-A\delta$.  The analogue of Theorem \ref{ethm15}
is:

\medskip
\begin{theorem}
\label{ethm17}
Denoting by $\rho^{(\delta,-)}$ the ``stationary profile'' when $R=1-A\delta$,  then for all $r\in [0,1)$
    \begin{equation}
     \label{e4.21}
\lim_{\delta\to 0} \rho^{(\delta,-)} (r)=2j(1-r) + \rho(1),\quad \rho(1):= \frac{j}{A}
   \end{equation}

\end{theorem}

\medskip

\noindent
{\bf Proof.} The main difference with  Theorem \ref{ethm15} is that now we have to deal
with an interval $[R,1]$ which
depends on $\delta$ and which shrinks to zero as $\delta\to 0$.  We can however set the problem
in such a way that the interval is the whole $[0,1]$ for all $\delta$.
To this end we introduce another map $T_\delta^{(\delta,-)}$  which, for a special choice of the parameters, will have the same fixed points as  $S_\delta^{(\delta,-)}$.  Given a non negative function $v$ we set
\begin{equation}
\label{e4.00}
T^{(\delta,-)}_\delta(u)= j\delta D_0 - v+ %\mathbf 1_{r\in [0,x]}
G_{{\delta}}^{\rm neum}*u
\end{equation}
Then a fixed point $u$ must have the form: $u=j\delta D_0- v +\psi(r) $ with $\psi$
such that
    \begin{equation}
     \label{e4.22}
\psi = G^{\rm neum}_\delta * [j\delta D_0- v + \psi ]
   \end{equation}
where, as before, $G^{\rm neum}_t$ is the Green function of the linear heat equation in $[0,1]$ with Neumann boundary conditions.

It is readily seen that if $\rho = \rho^{(\delta,-)}= \rho^{(\delta,-)}\mathbf 1_{r\notin[1-A\delta,1]}$ and
    \begin{equation}
     \label{e4.23}
v(r) := \mathbf 1_{r\in [1-A\delta,1]} G^{\rm neum}_\delta * [j\delta D_0+ \rho](r)
   \end{equation}
then
    \begin{equation}
     \label{e4.24}
\psi(r) := \begin{cases} \rho(r) & \text{if $r \in [0,1-A\delta]$}\\
v(r) & \text{if $r \in [1-A\delta,1]$}\end{cases}
   \end{equation}
solves \eqref{e4.22}.

On the other hand \eqref{e4.22} can be solved by iteration getting, analogously to \eqref{e4.4},
\begin{equation}
\label{e4.25}
\psi(r)=\sum_{n\ge 0} \{ j\delta G^{\rm neum}_{(n+1)\delta}(r,0) - \int_{1-A\delta}^1 G^{\rm neum}_{(n+1)\delta}(r,r') v(r') \}
\end{equation}
but again we need a proof that the series is convergent.  The Green function converges exponentially:
 \begin{equation}
\label{e4.26}
| G^{\rm neum}_{t}(r,r') - 1| \le c e^{-bt},\quad c>0,\;\;b>0
\end{equation}
Moreover, by its definition, see \eqref{e4.2.1},
\begin{equation}
\label{e4.27}
  \int_{1-A\delta}^1   v(r)dr = j\delta
\end{equation}
Then
\begin{equation*}
%\label{e4.25}
 \Big |j\delta G^{\rm neum}_{n\delta}(r,0) - \int_{1-A\delta}^1 G^{\rm neum}_{n\delta}(r,r') v(r') \Big |
 \le  c' e^{-bn \delta}
\end{equation*}
so that the series \eqref{e4.25} converges exponentially uniformly in $\delta$.
\vskip.3cm
Let us now add a superscript $(\delta,-)$ to $\psi$ and $v$ to underline their dependence on $\delta$.  We shall first prove that $\psi^{(\delta,-)}$ is equicontinuous:

\medskip
\begin{lemma}
For any $\eps>0$ there is $\alpha>0$ so that
for all $\delta$
\begin{equation}
\label{e4.28}
 \sup_{|r-r'|\le \alpha} |\psi^{(\delta,-)}(r)-\psi^{(\delta,-)}(r')| \le \eps
\end{equation}

\end{lemma}

\medskip

\noindent
{\bf Proof.}  By \eqref{e4.26} given any $\eps>0$ there is $T>0$ so that
\begin{equation}
\label{e4.29}
\sum_{n:n\delta \ge T} | j\delta G^{\rm neum}_{(n+1)\delta}(r,0) - \int_{1-A\delta}^1 G^{\rm neum}_{(n+1)\delta}(r,r'') v^{(\delta,-)}(r'') | \le \eps
\end{equation}
It is well known that for any $\zeta>0$ and $\tau>0$ there is $\alpha>0$ so that
\begin{equation*}
%\label{e4.28}
\sup_{t\ge \tau} \sup_{|r-r'|\le \alpha} \sup_{r''}|G^{\rm neum}_{t}(r,r'')-G^{\rm neum}_{t}(r',r'')| \le \zeta
\end{equation*}
By bounding $\dis{G^{\rm neum}_{t} (r,r')\le \frac{c}{\sqrt t}}$
\begin{eqnarray*}
|\psi^{(\delta,-)}(r)-\psi^{(\delta,-)}(r')| &\le &  2\eps + 4 \sum_{n\delta \le \tau} \frac{j\delta c}{\sqrt{ n\delta}}
+ \delta^{-1} T 2j\delta \zeta
\end{eqnarray*}
By choosing $\zeta = \eps/T$ and $\tau = \eps^2$ we then have  the right hand side  bounded proportionally to $\eps$ and the lemma is proved.  \qed

\vskip.5cm
By \eqref{e4.29} and the lemma we have that for any $\eps$ and for all $\delta$ small enough:
\begin{equation}
\label{e4.30}
 \Big|\psi^{(\delta,-)}(r) - j\delta \sum_n \Big(  G^{\rm neum}_{(n+1)\delta}(r,0) -
 G^{\rm neum}_{(n+1)\delta}(r,1)\Big ) \Big| \le \eps
\end{equation}
so that
\begin{equation}
\label{e4.31}
 \lim_{\delta\to 0} \psi^{(\delta,-)}(r) = \psi(r)= \int_0^\infty \{G^{\rm neum}_{t}(r,0) -
 G^{\rm neum}_{t}(r,1)\}
\end{equation}
which proves that $\rho^{(\delta,-)}(r)$ converges to $\psi(r)$ for all $r<1$.   As in the previous case
with $R<1$ fixed, the right hand side is identified to be a weak solution of the equation
\begin{equation}
\label{e4.33}
\psi'' + jD_0 -jD_1=0
\end{equation}
on $\mathbb R$ symmetric under all reflections of $[0,1]$.  To determine the solution we need
another condition, we are going to prove that at the right endpoint
\begin{equation}
\label{e4.32}
A \psi(1) = j
\end{equation}
Indeed,
\[
j\delta=  \int_{1-A\delta}^1   v^{(\delta,-)}(r)dr = A\delta v^{(\delta,-)}(1) + \int_{1-A\delta}^1   [v^{(\delta,-)}(r)-v^{(\delta,-)}(1)]dr
\]
Recalling \eqref{e4.24}, $v^{(\delta,-)}(r)=\psi^{(\delta,-)}(r)$, $r\in (1-A\delta,1)$, hence
\[
|j\delta -A\delta \psi^{(\delta,-)}(1)| \le  A\delta \sup_{1-A\delta\le r \le 1}|  \psi^{(\delta,-)}(r)-\psi^{(\delta,-)}(1)|
\]
By \eqref{e4.28} in the limit as $\delta\to 0$ we then obtain \eqref{e4.32}.  The weak solution of \eqref{e4.33} with the condition \eqref{e4.32} is the function on the right hand side of \eqref{e4.21}.   \qed

\vskip2cm

\subsection{Stationarity of the linear profiles}

\noindent
{\bf Proof of Theorem \ref{thme4}.}  To underline the choice of the initial datum we denote
the limit profile $\rho_t$ of Theorem \ref{Teo:Hydro} by
$\rho_t= S_t(\rho_{\rm init})$.  We fix $\tau>0$ and have by Theorem \ref{Teo:Char}
\begin{equation}
\label{e4.32aaa.1}
\lim_{n\to \infty} \Big| F(r;S_\tau(\rho^{(M)})) -F(r;S^{(\tau2^{-n},{-})}_\tau(\rho^{(M)}))\Big| =0,\quad \text{for all $r\in [0,1]$}
\end{equation}
Denote by $\rho^{(\tau2^{-n},-)}$ the stationary profile for the evolution $S_t^{(\tau2^{-n},-)}$
which converges to $\rho^{(M)}$, then
\begin{equation}
\label{e4.32aaa.2}
 \Big| F(r;S^{(\tau2^{-n},{-})}_\tau(\rho^{(\tau2^{-n},-)})) -F(r;S^{(\tau2^{-n},{-})}_\tau(\rho^{(M)}))\Big| \le
 \|S^{(\tau2^{-n},{-})}_\tau(\rho^{(\tau2^{-n},-)}))-S^{(\tau2^{-n},{-})}_\tau(\rho^{(M)}))\|_1
\end{equation}
Since $G^{\rm neum}_t$ is a contraction in $L_1$ as well as $K^{(\delta)}$ (see {Lemma (7.3) in \cite{CDGP1}})
we have
\begin{equation}
\label{e4.32aaa.3}
 \|S^{(\tau2^{-n},{-})}_\tau(\rho^{(\tau2^{-n},-)}))-S^{(\tau2^{-n},{-})}_\tau(\rho^{(M)}))\|_1
 \le \|\rho^{(\tau2^{-n},-)}-\rho^{(M)}\|_1
\end{equation}
which, {from Theorems \ref{ethm15} and \ref{ethm17},} vanishes  as $n\to \infty$.  Since $S^{(\tau2^{-n},{-})}_\tau(\rho^{(\tau2^{-n},-)})=\rho^{(\tau2^{-n},-)}$,
\begin{equation}
\label{e4.32aaa.4}
\lim_{n\to \infty} \Big| F(r;S_\tau(\rho^{(M)})) -F(r;\rho^{(\tau2^{-n},-)}))\Big| =0,\quad \text{for all $r\in [0,1]$}
\end{equation}
which concludes the proof because, as already observed,
\begin{equation}
\label{e4.32aaaAA.4}
\lim_{n\to \infty} \Big| F(r;\rho^{(M)}) -F(r;\rho^{(\tau2^{-n},-)})\Big| =0,\quad \text{for all $r\in [0,1]$}
\end{equation}
\qed

\section{Super-hydrodynamic limit}

The main result in this section is a proof of a  loss of memory of the initial conditions on
long hydrodynamic times.  This result will be obtained by introducing couplings and to this end it will be convenient to
label the particles.  We shall then conclude the section by using the loss of memory result
to prove convergence to linear stationary profiles and control the
super-hydrodynamic limit.

\vskip.5cm
\begin{definition} [Labeled configurations]
\label{defin:k5.1}

A labeled configuration is a pair $(\und x, I)$ where $I$ is a finite subset of $\mathbb N$
and $\und x$ a map from $I$ to $[0,\eps^{-1}]$: $I$ are the labels and $\und x$ the positions of the labeled particles.  We shall also write $\und x=\{x_i,\, i \in I\}$.
To any labeled configuration $(\und x,I)$ we associate the unlabeled
configuration $\xi_{\und x,I}$:
%\in \mathcal X$
%given by
	\begin{equation}
\label{1.11111111111}
\xi_{\und x,I}(x)=\sum_{i\in I} \mathbf 1_{x_i=x}
	\end{equation}
\end{definition}

\medskip

We shall couple the evolution starting from $(\und x_0,I_0)$ and $(\und y_0,J_0)$ where
$I_{ 0}=\{1,..,n\}$ and $J_{ 0}=(1,...,n+m)$, $n>0$, $m\ge 0$.  The coupled process will be
a jump Markov process %$(\und x_t, I_t, \und y_t,J_t,N_t)$, $t\ge 0$,
on a state space $S$ which is the family of all
$(\und x, I, \und y,J,N)$
such that
$I \subset J$, $J\setminus I$ has cardinality ${ \le}m$ and $N {=} \max \{i\in J\}$.

The coupled process starts from  $(\und x_0, I_0, \und y_0,J_0,n+m)$ and it is completely defined
once we specify the possible jumps and their intensities starting from  any element
$(\und x, I, \und y,J,N)$ in the state space  $S$. To this end we introduce the set
\[
I_{=}= \{ i\in I: x_i=y_i\}
\]
and call $(\und x', I', \und y',J',N')$ the elements after the jump.
The jumps are of four types:

\medskip

{\begin{itemize}
\item   {\em Single random walk jumps.} They are independent random walk jumps involving the restricted configurations $(\und x, I \setminus I_= )$ and $(\und y, J \setminus I_=)$. For
any of these jumps it will be $I'=I$, $J'=J$, $N'=N$.
The jumps indexed by $i\in I\setminus I_{=}$ are such that  $\und y'=\und y$ and  $x'_j=x_j$ for $ j\ne i$,
while $x_i \to x_i \pm 1$
with intensity $1/2$ and $x'_i=x_i\pm 1$ if this is
in $[0,\eps^{-1}]$, otherwise $x'_i=x_i$.
Analogously the jumps indexed by $i\in J\setminus I_{=}$ are such that $\und x'=\und x$ and
$y'_j=y_j$ for $ j\ne i$,   while
$y_i \to y_i \pm 1$ with intensity $1/2$
 and $y'_i=y_i\pm 1$ if this is in $[0,\eps^{-1}]$, otherwise $y'_i=y_i$.
We denote by $\mathcal L_s$ the Markov generator describing the single random walk jumps. It is given by:
\begin{eqnarray}
\mathcal L_s f(\und x,\und y) & = & \sum_{i \in I \setminus I_{=}}\frac 1 2 \left\{ \left( f(\und x^{i,+},\und y)-f(\und x, \und y)\right) \mathbf 1_{\{x_i \neq \eps^{-1}\}}+
\left( f(\und x^{i,-},\und y)-f(\und x, \und y)\right) \mathbf 1_{\{x_i \neq 0\}}\right\} \nn\\
& + & \sum_{i \in J \setminus I_{=}}\frac 12  \left\{\left( f(\und x,\und y^{i,+})-f(\und x, \und y)\right) \mathbf 1_{\{y_i \neq \eps^{-1}\}}+
\left( f(\und x,\und y^{i,-})-f(\und x, \und y)\right) \mathbf 1_{\{y_i \neq 0\}}\right\} \nn
\end{eqnarray}
where $\und x^{i,\pm}$ is the positions configuration obtained from $\und x$ by replacing $x_i$ with $x_i\pm 1$.  We have omitted to underline the dependence of $f$ on $I,J,N$ since they remain unchanged under the action of $\mathcal L_s$.
\item  {\em Double random walk jumps.} They are indexed by $i\in I_{=}$ and also for these jumps $I'=I$, $J'=J$, $N'=N$.  For each $i\in I_{=}$,  $x_i \to x'_i=x_i \pm 1$
and $y_i\to x'_i$, with intensity $1/2$ if $x_i\pm 1 \in [0,\eps^{-1}]$, otherwise the jump  is suppressed; all the other
positions are unchanged.
Let $\mathcal L_d$ be the Markov generator describing the double random walk jumps, then  it is given by:
\begin{eqnarray*}
&& \hskip-.5cm \mathcal L_d f(\und x,\und y)  = \nn \\ && = \sum_{i \in I_{=}} \frac 1 2\left\{\left( f(\und x^{i,+},\und y^{i,+})-f(\und x, \und y)\right) \mathbf 1_{\{x_i, y_i \neq \eps^{-1}\}}+
\left( f(\und x^{i,-},\und y^{i,-})-f(\und x, \und y)\right) \mathbf 1_{\{x_i,y_i \neq 0\}} \right\}
\end{eqnarray*}
\item  {\em Creation events.}  At rate $\eps j$, $N'= N+1$, $I'= I\cup\{N+1\}$, $J'= J\cup\{N+1\}$,
$x'_i=x_i, i\in I$, $x'_{N+1}=0$;  $y'_i=y_i, i\in  J$, $y'_{N+1}=0$. We call $\mathcal L_{\rm cr}$ the Markov generator associated to these events.
%\begin{eqnarray*}
%&&\hskip-.5cm\mathcal L_{\rm cr} f(\und x, I, \und y, J, N)  =   \nn \\
%&&=\eps  j \, \{f((\und x,0), I \cup \{N+1\}, (\und y,0), J \cup \{N+1\}, N+1)-f(\und x, I, \und y, J, N)\}
%\end{eqnarray*}
\item  {\em Death events.} At rate $\eps j$ both $I$ and $J$ loose an element while $N$ is unchanged.  The
configuration after the death event is obtained in two steps. In the first step
we erase from $\und x$
and $\und y$ their rightmost particle with largest label, say $x_i$  and  $y_j$. That is also the final step
 if $j\notin I$ or if $i=j$.
If instead $i\ne j$ and $j\in I$ we have two subcases: if $x_j \le y_i$ we relabel
$y_i$ as $y_j$ so that the label $i$ disappears from $I$ and $J$. If instead
$y_i < x_j$ we relabel
$x_j$ as $x_i$ so that the label $j$ disappears from $I$ and $J$. We denote by $\mathcal L_{\rm ann}$ the Markov generator associated to the death events.
\end{itemize}}

\medskip
\noindent
It directly follows from the above rules that:

\medskip

\begin{lemma}
\label{lemma4.1}
{In all the above cases $(\und x', I', \und y',J',N')\in S$ and
the set $I\setminus I_{=}$ does not increase after any of the above jumps.
Moreover in the case of a death event, if $i\in I \cap I'$,
 the interval with endpoints $x_i$ and $y_i$ may only change in such a way that the distance $|x_i-y_i|$ decreases.}
\end{lemma}

\medskip

One can then easily check that

\medskip

\begin{lemma}
\label{lemma4.2}
The above rules can be used to define a jump process with
state space $S$, denoted by
$(\und x(t), I(t), \und y(t),J(t),N(t))$. Its generator $\mathcal L$ is
	\begin{equation}
\label{n4.2}
\mathcal L = \mathcal L_s + \mathcal L_d + \mathcal L_{\rm cr} + \mathcal L_{\rm ann}
	\end{equation}
	%with $\mathcal L_s, \mathcal L_d, \mathcal L_{\rm cr}$ and $\mathcal L_{\rm ann}$  defined above.
where $\mathcal L_s$ describes the single random walk jumps; $\mathcal L_d$ the double random walk jumps; $ \mathcal L_{\rm cr}$ the creation  and $\mathcal L_{\rm ann}$ the annihilation jumps.
%
%
%\end{definition}
%
%

The processes
$\xi_{\und x(t),I(t)}$ and  $\xi_{\und y(t),J(t)}$ are then both  Markov  with generator
$L$ {defined in \eqref{generatore}}.

\end{lemma}

We say that $i$ is a discrepancy at time $t$ if it belongs to the set
     \begin{equation}
		\label{e4.3.1a}
 D_{\ne}(t) = I(t) \setminus I_{=}(t)=\Big\{ i \in I(t): x_i(t)\ne y_i(t)\Big\}
 %,\quad {\cal Z}(t)= |D_{\ne} (V_t(\und x,\und y,\om^*))|
      \end{equation}
By Lemma \ref{lemma4.1}, $D_{\ne}(t) \subset D_{\ne}(0){\subseteq I_0}$ hence, if $i\in D_{\ne}(t)$, then $i\in \{1,..,n\}$.
We denote by $|D_{\ne}(t) |$ the cardinality of
$D_{\ne}(t)$ which thus
counts the number of discrepancies at time $t$.

%
%Let also $V^0_t(\und x,\und y,\om^*)$ the
%evolution obtained without $\om_0$, so that
%$V^0_t(\und x,\und y,\om^*) =V_t(\und x,\und y,\om^*)$ for $t<t_{1;0}$.
%Thus $V^0_t$ describes the
%evolution of $n$ independent random walks
%starting from $\und x$ and $n$ independent random walks
%starting from $\und y$ with the rule that when
%particles with same label meet they stick together then after.
%Observe that

\medskip

\begin{lemma}
\label{zbound}
With the above notation, {for any $t \ge 0$ we have}
%riting $(\und x_t,\und y_t):= V_t(\und x,\und y,\om^*)$ and recalling {\eqref{1.11111111111}}
     \begin{equation}
		\label{e4.3.1a.1}
{\sum_{x=0}^{\eps^{-1}}} |\xi_{\und x(t),I(t)}(x)- \xi_{\und y(t),J(t)}(x)|  \le |D_{\ne}(t) | +m
      \end{equation}

\end{lemma}

\medskip

\noindent
{\bf Proof.}  Shorthand $\xi_t(x)=\xi_{\und x(t),I(t)}(x)$ and $\xi'_t(x)=\xi_{\und y(t),J(t)}(x)$. Then
	\begin{eqnarray*}
&&\hskip-1.5cm
\sum_{x=0}^{\eps^{-1}} |\xi_{t}(x)-\xi'_{t}(x)|=
\sum_{x=0}^{\eps^{-1}} \big |\sum_{i\in I(t)} \mathbf 1_{x_i(t)=x}- \sum_{i\in J(t)}\mathbf 1_{y_i(t)=x}\big|
\\&& \hskip1cm
\le \sum_{x=0}^{\eps^{-1}}\Big\{\sum_{i\in  D_{\ne}(t)}\big | \mathbf 1_{x_i(t)=x}- \mathbf 1_{y_i(t)=x}\big|
+\sum_{i\in J(t)\setminus I(t)} \mathbf 1_{y_i(t)=x}\Big\} \nn \\
&&\hskip3cm {=|D_{\ne}(t)| + |J(t)\setminus I(t)|}
	\end{eqnarray*}
{	then the result follows since $J(t)\setminus I(t) \subseteq J_0 \setminus I_0$ and $| J_0 \setminus I_0|=m$.}
\qed

\vskip.5cm
\noindent
Call $(\und x^0(t),\und y^0(t))$ the  independent random walk process starting from
$\und x^0(0)=(x_1,..,x_n)$ and $\und y^0(0)=(y_1,..,y_{n})$.
Call $\tau^0_i$, $i=1,..,n$, the first time $t$ when $x^0_i(t)=y^0_i(t)$ and
%
%when all the death and creation events
%are suppressed so that each pair $(x^0_i(t),y^0_i(t))$, $i=1,..,n$, moves independently of the others
%and the elements in each pair are independent random walks till when they meet, after that they stick together.
%The variables $y^0_i(t)$, $i=n+1,..,n+m$, are random walks independent of all the others.  We call
     \begin{equation}
		\label{n4.4}
 D^0_{\ne}(t) = \{ i \in \{1,..,n\}: \tau^0_i>t\}
% x^0_i(t)\ne y^0_i(t)\Big\}
%,\quad {\cal Z}(t)= |D_{\ne} (V_t(\und x,\und y,\om^*))|
      \end{equation}
and shall prove below that $| D^0_{\ne}(t)|$ stochastically bounds $| D_{\ne}(t)|$.

With this aim we introduce a process $(\und x(t), I(t), \und y(t),J(t),N(t);\und x^0(t),\und y^0(t))$
which couples the two processes $(\und x(t), I(t), \und y(t),J(t),N(t))$ and $(\und x^0(t),\und y^0(t))$.
We denote its generator by
           \begin{equation}
		\label{n4.6}
{\hat {\cal L} = \hat {\cal L}_s + \hat {\cal L}_d + \hat {\cal L}_{\rm cr} + \hat {\cal L}_{\rm ann} + \hat {\cal L}^0}
      \end{equation}
{$\hat {\cal L}_d$, $\hat {\cal L}_{\rm cr}$ and $\hat {\cal L}_{\rm ann} $}  are the  same  as
$\mathcal L_d$, $\mathcal L_{\rm cr}$ and $\mathcal L_{\rm ann} $ leaving unchanged
$\und x^0$ and $\und y^0$.
Also { $\hat {\cal L}_s$}   describes the same jumps as
$\mathcal L_s$ but it also changes $\und x^0$ and $\und y^0$ {with  the following rules.   For any $i \in I_0 \setminus I_=(t)$, if} $x_i \to \min\{x_i+1,\eps^{-1}\}$,
then also $x^0_i \to \min\{x^0_i+1,\eps^{-1}\}$ and, if  $x_i \to \max\{x_i-1,0\}$,
then also $x^0_i \to \max\{x^0_i-1,0\}$ (analogous rule for the $y$-jumps).
{The generator {$\hat {\cal L}^0$}  takes into account the
independent jumps of $ x^0_i$ and $y^0_i$ relative to the labels $i \in I_0 \cap I_=(t)$  which are not  been taken into account by {$\hat {\cal L}_s$}. As before, for any $i \in I_0 \cap I_=(t)$, if $x_i, y_i \to \min\{x_i+1,\eps^{-1}\}$,
then also $x^0_i, y_i^0 \to \min\{x^0_i+1,\eps^{-1}\}$ and, if  $x_i, y_i \to \max\{x_i-1,0\}$,
then also $x^0_i, y_i^0 \to \max\{x^0_i-1,0\}$.}

%the   discrepancies
%in the independent
%random walkers process.
%
%
%    \begin{equation}
%		\label{e4.3.1b}
% {\cal{Z}}^0(t)= |D_{\ne} (V^0_t(\und x,\und y,\om^*))|
%      \end{equation}
%and  $\{x^0_i(t)\},\{y^0_i(t)\}$
%the positions of the particles in $V^0_t(\und x,\und y,\om^*)$,
%then:

\medskip
	\begin{lemma}
\label{lemmae3.9}
If $\hat P$ is the law of the above   process with generator {$\hat {\cal L}$}, then
	\begin{equation}
	\label{n4.7}
\hat
P\big[D_{\ne}(t)\subset D^0_{\ne}(t)\big]=1\;\;\; \text{for all $t\ge 0$}
	\end{equation}
	
	\end{lemma}

\noindent
{\bf Proof.} Let us consider $i\in \{1,..,n\}$ and suppose (for the sake of
definiteness) that initially
$x_i<y_i$ (recalling that $x^0_i=x_i$ and $y^0_i=y_i$). Call $\tau_i$
the first time $t$ when either $i$ leaves $I(t)$ or $i$ enters into $I_{=}(t)$
We claim that $x_i(t)=x^0_i(t)$ and $y_i(t) \le y^0_i(t)$ for $t<\tau_i$ and since this implies
\eqref{n4.7} the claim will prove the lemma.  Indeed the jumps described by {$\hat {\cal L}_s$}
preserve such a property and if {$\hat {\cal L}_{\rm ann}$} involves the label $i$ (in the case
we are considering  it will still be present
after the jump event)
then $x_i$ is unchanged
and $y_i$ may only stay or decrease.
 \qed

\bigskip
As a direct consequence we have

\medskip

\begin{theorem}
\label{thm4.5}
There are positive constants $c$ and $b$ so that for any {$t\ge 0$}, any $n$, $m$ and  any
initial configurations $\und x$ and $\und y$ as above
	\begin{equation}
	\label{3.2.001}
\hat E\big[|D_{\ne}(t)|\big]  \le cn e^{- b \eps^2 t}
	\end{equation}
($\hat E$ denoting expectation with respect to the measure $\hat P$).

\end{theorem}

\medskip

\noindent
{\bf Proof.}  	By \eqref{n4.7} it is enough to prove the inequality
for $E^0\big[|D^0_{\ne}(t)|\big]$, $E^0$ the expectation for the   independent walkers process.
The bound will follow from the inequality
\begin{equation}\label{AAA}
{p_t {=p_t(i)}:=}P^0[\tau^0_i>t ] \le c e^{- b \eps^2 t}
\end{equation}
for any $i \in I_0=\{1, \dots, n\}$. There is $\ga>0$ so that
supposing $x_i< y_i$
   $$
p_{ \eps^{-2}}\ge P^0\Big[x^0( \eps^{-2})\ge \frac {y+x}2, \: \:
y^0( \eps^{-2})\le \frac {y+x}2\Big]=\Big(P^0\Big[x^0( \eps^{-2})\ge \frac {y+x}2\Big]
\Big)^2 \ge \ga
   $$
hence
\begin{equation}
\label{3.2.001.1}
p_{t}\le  (1- \ga)^{\eps^2 t -1} = c e^{-b\eps^2 t},\qquad b =- \log (1-\ga),\qquad c=(1-\ga)^{-1}
\end{equation}
{Now we have
\begin{eqnarray}
 P^0\big[|D_{\ne}(t)|=k\big] = \sum_{{\cal I}\subseteq I_0: \, |{\cal I}|=k}  \quad  \prod_{ i\in {\cal I}}  P^0\big[\tau^0_{i}>t\big] \cdot \prod_{j\notin {\cal I}}  P^0\big[\tau^0_{j}\le t \big] \nn \\
= \binom{n}{k} \, p_t^k \, (1-p_t)^{n-k}
\end{eqnarray}
then  $E^0\big[|D_{\ne}(t)|\big] = n \, p_t$, this proves the Theorem.}
\qed

\vskip.5cm

\subsection{Convergence to linear profiles}
\label{undici}

\vskip.3cm

We start by proving  Theorem \ref{thme2.3}, to this end we
show that two initial profiles
with the same mass  (or two initial configurations with the same total
number of particles) become indistinguishable on the hydrodynamic
time scale. %This is contained in the propositions which follow below and
%which are the main ingredients for the proof of Theorem \ref{thme2.3}.

\vskip.5cm

\begin{proposition}[Loss of memory for $\rho_t$]
\label{prop-squeezing}
Let  $\rho_{\rm init},\tilde{\rho}_{\rm init}$ be as in Definition \ref{def:2.1}.  Suppose  $F(0;\rho_{\rm init})=F(0;\tilde\rho_{\rm init})=:M$,
 then
 %for any $r \in [0,1]$,
  \begin{equation}
  \label{loss}
\lim_{t \to \infty} \sup_{r \in [0,1]}\big|F(r; S_t (\rho_{\rm init}))-F(r; S_t (\tilde \rho_{\rm init}))\big| = 0
  \end{equation}

\end{proposition}
\medskip
\noindent
{\bf Proof.}  We shall use a corollary of Theorem \ref{Teo:Hydro} which may have
an interest in its own right.  Let  $\rho_{\rm init}$, $\xi$  and  $\rho_t(r)$ as in
Theorem \ref{Teo:Hydro} then for any $t>0$
   \begin{equation}
   \label{5.15.00}
\lim_{\eps \to 0}E^{(\eps)}_{\xi}\Big[\max_{x\in [0,\eps^{-1}]}|\eps F_\eps(x;\xi_{\eps^{-2} t}) - F(\eps x;\rho_{t})| \Big] = 0
    \end{equation}
{\em Proof of \eqref{5.15.00}.} For any $\zeta>0$  define
\begin{eqnarray*}
&&\mathcal E_\zeta(\eps,t):= \Big\{{\max_{x \in [0,\eps^{-1}]}}\big|F(\eps x; S_t(\rho_{\rm init}))- \eps F_\eps({ x}; \xi_{\eps^{-2}t})\big|\le \zeta \Big\} %,\;\;  x_r:=[\eps^{-1} r]
\end{eqnarray*}
Then, {from the Cauchy-Schwarz inequality we have}
   \begin{eqnarray*}
%   \label{5.15.00}
%\lim_{\eps \to 0}
&& E^{(\eps)}_{\xi}\Big[\max_{x\in [0,\eps^{-1}]}|\eps F_\eps(x;\xi_{\eps^{-2} t}) - F(\eps x;\rho_{t})| \Big] \nn \\
&&{ \le   \zeta P^{(\eps)}_{\xi}\Big[\mathcal E_\zeta(\eps,t)\Big] + P^{(\eps)}_{\xi}\Big[\mathcal E_\zeta(\eps,t)^c\Big]^{1/2}E^{(\eps)}_{\xi}\Big[\Big(\max_{x\in [0,\eps^{-1}]}|\eps F_\eps(x;\xi_{\eps^{-2} t}) - F(\eps x;\rho_{t})|\Big)^2\Big]^{1/2} }\nn \\
&& \hskip5.5cm\le \zeta + P^{(\eps)}_{\xi}\Big[\mathcal E_\zeta(\eps,t)^c\Big]^{1/2}E^{(\eps)}_{\xi}\Big[(\eps |\xi_{\eps^{-2} t}|+M)^2\Big]^{1/2}
    \end{eqnarray*}
{because $F(0,\rho_t)=F(0,\rho_{\rm init})=M$.}
By  Theorem \ref{Teo:Hydro}  $P^{(\eps)}_{\xi}\Big[\mathcal E_\zeta(\eps,t)^c\Big]$ vanishes while
by Theorem \ref{thme3.6a}, $E^{(\eps)}_{\xi}\Big[(M+\eps |\xi_{\eps^{-2} t}|)^2\Big]\le c$ uniformly in $\eps$.
\eqref{5.15.00} is thus proved.

\medskip
{Let $\{\tilde \xi\}$ be the family of initial data which approximate $\tilde \rho_{\rm init}$,
chosen in such a way that for all $\eps$, $|\tilde \xi|=| \xi|= :n_\eps$. Calling $x_r:=[\eps^{-1} r]$,
$r \in [0,1]$, since $S_t(\rho_{\rm init}), S_t(\tilde \rho_{\rm init})$ are bounded we have
  \begin{eqnarray*}
%  \label{loss}
&&\big|F(r; S_t (\rho_{\rm init}))-F(r; S_t (\tilde \rho_{\rm init})) \big| \le  c \, |r-\eps x_r |+E^{(\eps)}_{\xi}\Big[|F(\eps x;S_t(\rho_{\rm init}))-\eps F_\eps(x;\xi_{\eps^{-2} t})| \Big]\nn \\
&& +E^{(\eps)}_{\tilde \xi}\Big[|F(\eps x;S_t(\tilde\rho_{\rm init}))-\eps F_\eps(x;\tilde\xi_{\eps^{-2} t})| \Big]+\Big|
E^{(\eps)}_{\xi}\Big[\eps F_\eps(x_r;\xi_{\eps^{-2} t})\Big] - E^{(\eps)}_{{ \tilde \xi}}\Big[\eps F_\eps(x_r;\tilde \xi_{\eps^{-2} t})\Big]\Big|
  \end{eqnarray*}
for some $c\ge 0$, then,  by \eqref{5.15.00},
\eqref{e4.3.1a.1} with $m=0$ and
\eqref{3.2.001}
  \begin{eqnarray}
  \label{lossbis}
\Big|F(r; S_t (\rho_{\rm init}))-F(r; S_t (\tilde \rho_{\rm init}))\Big| \le \lim_{\eps \to 0} \Big|
E^{(\eps)}_{\xi}\Big[\eps F_\eps(x_r;\xi_{\eps^{-2} t})\Big] - E^{(\eps)}_{{ \tilde \xi}}\Big[\eps F_\eps(x_r;\tilde \xi_{\eps^{-2} t})\Big]\Big| &&\nn \\
\le  \lim_{\eps\to 0} c\eps n_\eps e^{- b  t} \le c F(0;\rho_{\rm init}) e^{- b  t} &&
  \end{eqnarray}
\eqref{loss} is then proved. }\qed

\vskip1cm

\noindent
{\bf Proof of Theorem \ref{thme2.3}.}
Equation \eqref{e2.666}  follows from  Theorem \ref{thme4} and  Proposition \ref{prop-squeezing} with $\tilde \rho_{\rm init}=\rho^{(M)}$.  \qed

\vskip.5cm

We fix arbitrary $M>0$ as an upper bound for the total {macroscopic} mass with $\eps^{-1}M$
bounding the total number of particles.

\medskip

\begin{definition}

For any $\eps>0$ and any positive
integer
$N\le M\eps^{-1}$ we denote by
%
%We define a family of configurations $\{\eta^{(N,\eps)}\}$, $N\le M\eps^{-1}$, $M>0$,
$\eta^{(N,\eps)} \in \mathbb N^{\Lambda_\eps}$ the following particle approximation of the invariant profile
$\rho^{(\eps N)}$.  We set  $\eta^{(N,\eps)}(\eps^{-1})=0$ and define iteratively for any $x\in [0,\eps^{-1}-1]$:
\begin{eqnarray}
\label{54}
\sum_{y=0}^x\eta^{(N,\eps)}(y)= {\bigg\lceil\eps^{-1}\int_0^{\eps(x+1)}\rho^{(\eps N)} \bigg\rceil}
\end{eqnarray}
{where $\lceil z \rceil$ is the  smallest integer $\ge z$.}

\end{definition}

Observe that $\dis{\sum_{y=0}^{\eps^{-1}-1}\eta^{(N,\eps)}(y) = N}$ and that for any $m>0$   the sequence  $\eta^{([\eps^{-1}m],\eps)}$ satisfies the conditions
in Definition \ref{def:2.1} with respect to $\rho_{\rm init}= \rho^{(m)}$.

\medskip

\begin{proposition}

For any $\zeta>0$ and $M>0$ there are $t$ and $\eps^*$ so that for any $\eps\le \eps^*$:
\begin{eqnarray}
\label{55}
\sup_{\xi: |\xi|\le M\eps^{-1}} E^{(\eps)}_\xi
\Big[ \max_{x {\in [0,\eps^{-1}]}} \big|\eps F_\eps(x;\xi_{\eps^{-2}t} )
- F(\eps x; \rho^{(\eps |\xi|)})\big|\Big] \le \zeta
\end{eqnarray}

\end{proposition}

\medskip
\noindent
{\bf Proof.}  We split the interval $[0,M]$ into intervals of length $\theta$,
$\theta>0$, calling $\theta_n=n\theta$. We choose $\theta$ so small that
\[
\max_{n} \sup_{m\in [\theta_n,\theta_{n+1}]}
\int_0^1 |  \rho^{( \theta_n)} (r)-  \rho^{(m)} (r)|
\le  {\frac \zeta 2}
\]
{Let $\eta^{([\eps^{-1}\theta_n],\eps)}_t$ be the process with generator \eqref{generatore} and initial configuration $\eta^{([\eps^{-1}\theta_n],\eps)}$, then $|\eps|\eta^{([\eps^{-1}\theta_n],\eps)}|-\theta_n|\le \eps$.}  By \eqref{e4.3.1a.1} and Theorem \ref{thm4.5}   for any $n$ and any $\xi$ such that
 $\eps |\xi| \in [\theta_n,\theta_{n+1}]$,
\begin{eqnarray}
\label{57}
%\max_{x}
\hat E\Big[ \eps \sum_x \big|{\eta^{([\eps^{-1}\theta_n],\eps)}_{\eps^{-2}t} }(x)-
\xi_{\eps^{-2}t}(x)\big|\Big] \le \theta  {+\eps} +c M e^{-bt} \le \frac \zeta 4
\end{eqnarray}
%where $\hat E$ denotes expectation with respect to the coupling of
%the processes starting from $\eta^{([\eps^{-1}\theta_n],\eps)}$ and $
%\eta^{([\eps^{-1}\theta'],\eps)}$.
The last inequality requires $t$ large enough: $c M e^{-bt} <\zeta/8$ { and $\theta$ and $\eps$ small enough so that $\theta+\eps\le  \zeta/8$.}
By Theorem \ref{Teo:Hydro} and
\eqref{5.15.00},
since
$S_t(\rho^{(m)})=\rho^{(m)}$,
there is $\eps_1(t,\zeta;\theta)$ so that for all $\eps\le \eps_1(t,\zeta;\theta)$
\begin{eqnarray}
\label{56}
\max_{n} \quad
E^{(\eps)}_{\eta^{([\eps^{-1}\theta_n],\eps)}}\Big[ \max_{x {\in [0,\eps^{-1}]}} \big|\eps F_\eps(x;{\eta^{([\eps^{-1}\theta_n],\eps)}_{\eps^{-2}t}})
- F(\eps x; \rho^{( \theta_n)})\big|\Big] \le \frac \zeta 4
\end{eqnarray}
As a consequence
\begin{eqnarray*}
&& \hskip-2cm \sup_{\xi: \eps |\xi|\le M} \: E^{(\eps)}_\xi
\Big[ \max_{x {\in [0,\eps^{-1}]}} \big|\eps F_\eps(x;\xi_{\eps^{-2}t} )
- F(\eps x; \rho^{(\eps |\xi|)})\big|\Big] \\&& \le
\sup_{\xi: \eps |\xi|\le M} E^{(\eps)}_\xi
\Big[ \max_{x {\in [0,\eps^{-1}]}} \big|\eps F_\eps(x; {\eta^{([\eps^{-1}\theta_n],\eps)}_{\eps^{-2}t}} )
- F(\eps x; \rho^{(\eps |\xi|)})\big|\Big] +\frac \zeta 4\\
%\hskip1cm && +
%\max_{n} \sup_{\xi: \eps |\xi|\in [\theta_n,\theta_{n+1}]}
%\hat E\Big[ \eps \sum_x|\eta^{([\eps^{-1}\theta_n],\eps)}_{\eps^{-2}t} (x)-
%\xi_{\eps^{-2}t}(x)|\Big] +
%\sup_{\xi: \eps |\xi|\le M} E^{(\eps)}_\xi
%\Big[ \max_{x} |\eps F_\eps(x; \eta^{([\eps^{-1}\theta_n],\eps)}_{\eps^{-2}t} )
%- F(\eps x; \rho^{(\eps |\xi|)})|\Big]
 \\
\hskip1cm && \le  \frac \zeta 2 + \max_{n} \sup_{m\in [\theta_n,\theta_{n+1}]}
\int_0^1 |  \rho^{( \theta_n)} (r)-  \rho^{(m)} (r)| <\zeta
\end{eqnarray*}
{this concludes the proof.}  \qed

\vskip.5cm

\subsection{Evolving profiles}

\vskip.4cm
{\bf Proof of Theorem \ref{thme3.6a}.}
{From the definition of the generator \eqref{generatore} we infer
that the induced process $|\xi_t|$ counting the number of particles at time $t$
evolves with the generator
\begin{equation}\label{Gen}
{\cal L}^{(\eps)}f(|\xi|) = j\eps  \{\Big( f(|\xi|+1)-f(|\xi|)\Big) + \mathbf 1_{|\xi|>0}
\Big( f(|\xi|-1)-f(|\xi|)\Big)\}
\end{equation}
acting on bounded functions $f:\mathbb N \to \mathbb R$.
Such generator is immediately recognized to be the generator
of the continuous time symmetric random walk on $\mathbb N$
at rate $j\eps$ and reflected at the origin.} \qed

\vskip.5cm

We also have that calling $\mathcal P^{(\eps)}_x$, $x\in \mathbb N$, the law of the random walk $x_t$ with generator $\mathcal L^{(\eps)}$
starting from $x$:
%
%Before  proceeding with the proof of the first statement of Theorem \ref{Teo:Local}, i.e. the super-hydrodynamical limit in \eqref{SuperHydro}, we need  to control the fluctuation of  $M_t^{(\eps)}$.
%We denote by $P_{s,\xi_*}^{(\eps)}$ the law of the process $\xi_t$ on $[0,\eps^{-1}]$, starting at time $s$ from the configuration $\xi_*$ (\colo{thus}  $P_{\xi}^{(\eps)}=P_{0,\xi}^{(\eps)}$).

\medskip

   \begin{lemma}
   \label{coro:B}

  Let $M'>0$ and $T>0$ then for any $\delta>0$ there is $M$ so that for all $\eps$ small enough, any
  $x \le \eps^{-1}M'$
and any $t \le \eps^{-3}T$,
  \begin{equation}
\label{ultimo}
\sup_{t\le \eps^{-2}T}\mathcal P^{(\eps)}_{x} \left[|x_t -x| \le \delta  \right] \ge 1-\delta,\qquad
 \sup_{t\le \eps^{-3}T}\mathcal P^{(\eps)}_{x} \left[x_t \ge \eps^{-1}M  \right] \le \delta
   \end{equation}
%
%   $x \le M\eps^{-1}$,
%
%  for any $\eps>0
%For any $\xi_*$ \colo{and for any} $s,T>0$ and $\alpha \in (0,1)$,
%  \begin{equation}
%\label{ultimo}
%\lim_{\eps \to 0}P^{(\eps)}_{\eps^{-3}s,\xi_*} \left[\big|M^{(\eps)}_{s+\eps T}-M^{(\eps)}_s\big|>\eps^{\frac{\alpha}2} \right] = 0
%   \end{equation}

   \end{lemma}

   %
%
%\medskip
%\noindent
%{\bf Proof.}
%It is sufficient to  prove the result for $s=0$ because $M^{(\eps)}_t$ is
%%an homogeneous Markov
%\colo{a process with stationary increments} , i.e. we prove that, for any initial configuration $\xi$,
%\begin{equation}
%\label{ultimo}
%\lim_{\eps \to 0}P^{(\eps)}_{\xi} \left[\big|M^{(\eps)}_{\eps T}-\eps \, |\xi|\big|>\eps^{\frac{\alpha}2} \right] = 0
%\end{equation}
%{Using $\frac{d}{dt}E_{\xi}[f(|\xi_t|)]=E_{\xi}[\mathcal Lf(|\xi_t|)]$ with $\mathcal L$ given by \eqref{Gen}
%  and  $f(n)=(n-|\xi|)^2$ we get}
%\begin{equation}\label{exp}
%\frac d{dt}\, E_{\xi}[(|\xi_t|-|\xi|)^2]=2\eps j -\eps j (1+ 2 |\xi|)P_\nu \left[|\xi_t|=0\right] \le 2j \eps
%\end{equation} Hence
%\begin{equation}\label{Var}
%E_{\xi}[(|\xi_t|-|\xi|)^2] \le 2\eps j t
%\end{equation}
%We recall that $M^{(\eps)}_{\eps T}=\eps |\xi_{\eps^{-2}T}|$. Therefore, by applying Markov inequality
%and using the bound \eqref{Var} one has
%\begin{equation}
%P^{(\eps)}_{\xi} \left[\big|M^{(\eps)}_{\eps T}-\eps \, |\xi|\big|>\eps^{\frac{\alpha}2} \right] =
%P^{(\eps)}_{\xi} \left[ (|\xi_{\eps^{-2}T}|-|\xi |)^2>\eps^{\alpha-2} \right] \le  2j\eps \,\frac{\eps^{-2} T}{\eps^{\alpha-2}} = c \eps^{1-\alpha}\;.
%\end{equation}
%\qed
%

\vskip1cm

\vskip.5cm
\noindent
{\bf Proof of Theorem \ref{Teo:Local}.}
The last statement of the theorem, i.e. that
\begin{equation}
M_t^{(\eps)}:= \eps \big|\xi_{\eps^{-3}t}\big| \to {B_{jt}} \qquad \text{as }\eps \to 0 \qquad \text{ \emph{in law}}
\end{equation}
with {$(B_t)_{t\ge 0}$}  the brownian motion on $\mathbb R_+$ with reflection at the origin, starting from ${B_0}= \lim_{\eps\to 0} M_0^{(\eps)}=\lim_{\eps \to 0}\eps |\xi|$,
follows from  Theorem \ref{thme3.6a}  and the fact that the diffusive scaling limit of the random walk is Brownian motion.

$\bullet$\; Subcritical regime.  \eqref{e2.6666xx} follows  directly from \eqref{55}.

$\bullet$\; Critical regime.  Let $t^*=\eps^{-3}t- s$, then by Lemma \ref{coro:B} for any given $s>0$,
with probability $\ge 1-\delta$, $|\xi_{t^*} \le \eps^{-1}M$.  By \eqref{55}, choosing $s$ large enough
 in the set $|\xi_{t^*} |\le \eps^{-1}M$,
\begin{eqnarray*}
\label{55}
 E^{(\eps)}_{\xi_{t^*}}
\Big[ \max_{x} |\eps F_\eps(x;\xi_{\eps^{-2}s} )
- F(\eps x; \rho^{(\eps |\xi_{t^*}|)})|\Big] \le \frac \zeta 2
\end{eqnarray*}
On the other hand by \eqref{ultimo} for $\eps$ small enough
  \begin{equation*}
%\label{ultimo}
P^{(\eps)}_{\xi_{t^*}} \left[||\xi_{t^*}| -|\xi_{t^*+\eps^{-2}s}| |\le \delta  \right] \ge 1-\delta
   \end{equation*}
so that \eqref{SuperHydro} follows from the continuity in $m$ of $F(0; \rho^{(m)})$.  \qed

\vspace{0.5cm}

{\bf Acknowledgments.} The research has been partially supported by PRIN 2009 (prot.
2009TA2595-002) and FIRB 2010 (grant n. RBFR10N90W). A. De Masi and E. Presutti
acknowledge kind hospitality at Universit\`a di Modena e Reggio Emilia.
G. Carinci and C. Giardin\`a  thanks Universit\`a dell'Aquila for welcoming during their visit.

\end{document}